\documentclass[a4paper,12pt,intlimits,oneside,reqno]{amsart}
\usepackage{enumerate}
\usepackage{amsfonts,amsmath, amsxtra,amssymb,latexsym, amscd,amsthm}

\usepackage{latexsym,amssymb}
\newtheorem{theorem}{Theorem}[section]
\newtheorem{lemma}[theorem]{Lemma}

\theoremstyle{corollary}

\theoremstyle{definition}
\newtheorem{definition}[theorem]{Definition}

\theoremstyle{remark}

\numberwithin{equation}{section}

\newcommand{\comment}[1]{}

\begin{document}

\title [Commutator of multilinear Hausdorff operator]{Weighted estimates for commutators of multilinear Hausdorff operators on variable exponent Morrey-Herz type spaces}

\author{Nguyen Minh Chuong}

\address{Institute of mathematics, Vietnamese  Academy of Science and Technology,  Hanoi, Vietnam.}
\email{nmchuong@math.ac.vn}

\author{Dao Van Duong}
\address{Shool of Mathematics,  Mientrung University of Civil Engineering, Phu Yen, Vietnam.}
\email{daovanduong@muce.edu.vn}

\author{Kieu Huu Dung}
\address{Shool of Mathematics, University of Transport and Communications- Campus in Ho Chi Minh City, Vietnam.}
\email{khdung@utc2.edu.vn}
\begin{abstract}
In this paper, we establish the boundedness of the commutators of multilinear Hausdorff operators on the product of some weighted Morrey-Herz type spaces with variable exponent with their symbols belong to both Lipschitz space and  central BMO space. By these, we generalize and strengthen some previous known results.
\end{abstract}

\maketitle

\section{Introduction}\label{section1}
Given $\Phi$ be a locally integrable function on $\mathbb R^n$. The $n$-dimensional Hausdorff operator $H_{\Phi,A}$ \cite{BM} is defined by
\begin{equation}\label{Hausdorff1}
H_{\Phi, A}(f)(x)=\int\limits_{\mathbb R^n}{\frac{\Phi(t)}{|t|^n}f(A(t) x)dt},\,x\in\mathbb R^n,
\end{equation}
where $A(t)$ is an $n\times n$ invertible matrix for almost everywhere $t$ in the support of $\Phi$.
It is well known that  if the function $\Phi$ and the matrix $A$ are taken appropriately, then the Hausdorff operator $H_{\Phi,A}$ reduces to many classcial operators in analysis, for example, the Hardy operator, the Ces\`{a}ro operator, the Hardy-Littlewood average operator and the Riemann-Liouville fractional integral operator. Some of their results have been significantly seen in \cite{BM}, \cite{CDH2016}, \cite{CHH2016}, \cite{FLL2009}, \cite{Moricz2005}, \cite{TXZ2011}, \cite{Xiao2001} and references therein.
\vskip 5pt
In addition, it is natural to extend the study on the linear operator to multilinear operator, which is actually necessary. Thus, the authors of this paper in \cite{HausdoffCDD} have recently investigated the multilinear operators of Hausdorff type $H_{\Phi,\vec A}$ given as follows:
\begin{equation}\label{mulHausdorff}
{H_{\Phi ,\vec{A} }}(\vec{f})(x) = \int\limits_{{\mathbb R^n}} {\frac{{\Phi (t)}}{{{{\left| t \right|}^n}}}} \prod\limits_{i = 1}^m {{f_i}} ({A_i}(t)x)dt,\,x\in\mathbb R^n,
\end{equation}
where $\Phi:\mathbb R^n \to [0,\infty)$ and $A_i(t)$ (for $i=1,..., m$) are $n\times n$  invertible matrices for almost everywhere $t$ in the support of $\Phi$, and $f_1, f_2, ..., f_m:\mathbb R^n\to \mathbb C$ are measurable functions and $\vec{f}=\left(f_1, ..., f_m\right)$ and $\vec{A}=\left(A_1, ..., A_m\right)$.  It is useful to remark that the weighted multilinear Hardy operators \cite{FGLY2015} and weighted multilinear Hardy-Ces\`{a}ro operators \cite{CHH2016} are two special cases of the multilinear Hausdorff operators $H_{\Phi,\vec A}$.
\vskip 5pt
\begin{definition}
Let $\Phi, \vec A, \vec f$ be as above. The Coifman-Rochberg-Weiss type commutator of multilinear Hausdorff operator is defined by
\begin{equation}\label{commulHausdorff}
{H_{\Phi ,\vec{A} }^{\vec b}}(\vec{f})(x) = \int\limits_{{\mathbb R^n}} {\frac{{\Phi (t)}}{{{{\left| t \right|}^n}}}} \prod\limits_{i = 1}^m {{f_i}} ({A_i}(t)x)\prod\limits_{i=1}^m{\big(b_i(x)-b_i(A_i(t)x)\big)}dt,\,x\in\mathbb R^n,
\end{equation}
where $\vec b= \left(b_1, ..., b_m\right)$ and $b_i$ are locally integrable functions on $\mathbb R^n$ for all $i=1,...,m$.
\end{definition}
Moreover, if we now take $m=n \geq 2, \Phi(t)=|t|^m.\omega(t)\chi_{[0,1]^m}(t)$ and $A_i(t)= t_i.I_m$ ($I_m$ is an identity matrix), for $t = (t_1, t_2, ..., t_m)$, where $\omega:[0,1]^m\to [0,\infty)$ is a measurable function, then $H_{\Phi,\vec A}^{\vec b}$ reduces to the  commutator of weighted multilinear Hardy operator due to  Fu et al. \cite{FGLY2015} defined as the following
\begin{equation}\label{Upsi}
{H_{\omega }^{\vec b}}(\vec{f})(x) = \int\limits_{{[0,1]^m}}\prod\limits_{i = 1}^m {{f_i}} ({t_i}x)\prod\limits_{i=1}^m{\big(b_i(x)-b_i(t_ix)\big)}\omega(t)dt,\,x\in\mathbb R^m.
\end{equation}
Also, by $\Phi(t)=|t|^n.\psi(t)\chi_{[0,1]^n}(t)$ and $A_i(t)= s_i(t).I_n$, where $\psi:[0,1]^n\to [0,\infty), s_1,...,s_m:[0,1]^n\to \mathbb R$ are  measurable functions, it is clear to see that $H_{\Phi,A}^{\vec b}$ reduces to the commutator of  multilinear Hardy-Ces\`{a}ro operator $U_{\psi,\vec {s}}^{m, n,\vec b}$ introduced by  Hung and Ky \cite{HK2015} as follows
\begin{equation}\label{MulUpsib}
U_{\psi,\vec {s}}^{m, n,\vec b}(x) = \int\limits_{{[0,1]^n}}\prod\limits_{i = 1}^m {{f_i}} ({s_i}(t)x)\prod\limits_{i=1}^m{\big(b_i(x)-b_i(s_i(t)x)\big)}\psi(t)dt,\,x\in\mathbb R^n.
\end{equation}
\vskip 5pt
In recent years, the theory of  function spaces with variable exponents has attracted much more the interest from many mathematicians (see, e.g., \cite{Almeida2012}, \cite{Bandaliev2010}, \cite{HausdoffCDD}, \cite{Diening1},  \cite{Guliyev}, \cite{I2010Hiro}, \cite{WZ2016} and others). It is interesting to see that this theory has had some important applications to the electronic fluid mechanics, elasticity, fluid dynamics, recovery of graphics, harmonic analysis and partial differential equations (see \cite{Almeida}, \cite{Chuong},  \cite{Chuong2016}, \cite{CUF2013}, \cite{Diening}, \cite{H2000}, \cite{Jacob}). 
\vskip 5pt
Let $b\in BMO(\mathbb R^n)$ and $T$ be a Calder\'{o}n-Zygmund singular integral operator with rough kernels.
From classical result of Coifman, Rochberg, and Weiss \cite{CRW1976},  Karlovich and Lerner \cite{KL2005} developed the boundedness of commutator $[b, T]$ to generalized $L^p$ spaces with variable exponent. Also, in order to generalize the result of  Chanillo \cite{C1982}, Izuki \cite{I2010} established the boundedness of the higher order commutator on Herz spaces with variable exponent.
\vskip 5pt
More recently, Wu \cite{W2014} considered the $m$th-order commutator for the fractional integral as follows
\[
I_{\beta, b}^m(f)(x)=\int\limits_{\mathbb R^n}{\frac{f(y)\big(b(x)-b(y)\big)^m}{|x-y|^{n-\beta}}dy},
\]
where $\beta\in (0,n)$, $b\in BMO(\mathbb R^n)$, $m\in\mathbb N$. Then  the author established the boundedness
for commutators of fractional integrals on Herz-Morrey spaces with variable exponent. 
\vskip 5pt
Motivated by above mentioned results,  the goal of this paper is to establish the boundedness for commutators of multilinear Hausdorff operators on the product of weighted Lebesgue, central Morrey, Herz, and Morrey-Herz spaces with variable exponent with their symbols belong to both Lipschitz spaces and central BMO spaces.
\vskip 5pt
Our paper is organized as follows. In Section 2, we give necessary preliminaries on weighted Lebesgue spaces, central Morrey spaces, Herz spaces, Morrey-Herz spaces with variable exponent and Lipschitz spaces, central BMO spaces. In Section 3, our main theorems are given. Finally, the results of this paper are proved in Section 4.
\section{Preliminaries}\label{section2}
In this section, let us recall some basic facts and notations which will be used throughout this paper.  The letter $C$ denotes a positive constant which is independent of the main parameters, but may be different from line to line. Given a measurable set $\Omega$, let us denote by $\chi_\Omega$ its characteristic function, by $|\Omega|$ its Lebesgue measure. For any $a\in\mathbb R^n$ and $r>0$, we denote by $B(a,r)$ the ball centered at $a$ with radius $r$.
\vskip 5pt
Next, we write $a \lesssim b$ to mean that there is a positive constant $C$, independent of the main parameters, such that $a\leq Cb$. Besides that, we denote $\chi_k=\chi_{C_k}$, $C_k=B_k\setminus B_{k-1}$ and $B_k = \big\{x\in \mathbb R^n: |x| \leq 2^k\big\}$, for all $k\in\mathbb Z$. 
\vskip 5pt
Now, we present the definition of the Lebesgue  space with variable exponent. For further readings on its deep applications in harmonic analysis, the interested reader may find in the works \cite{CUF2013}, \cite{Diening} and \cite{Diening1}.
\begin{definition}
Let $\mathcal P(\mathbb R^n)$ be the set of all measurable functions $p(\cdot):\mathbb R\to [1,\infty]$. For $p(\cdot)\in \mathcal P(\mathbb R^n)$, the variable exponent Lebesgue space $L^{p(\cdot)}(\mathbb R^n)$ is the set of all complex-valued measurable functions $f$ defined on $\mathbb R^n$ such that there exists constant $\eta >0$ satisfying
\[
F_{p}(f/\eta):=\int\limits_{\mathbb R^n\setminus\Omega_\infty}\left({\frac{|f(x)|}{\eta}}\right)^{p(x)}dx +\big\|f/\eta\big\|_{L^{\infty}(\Omega_\infty)}<\infty,
\]
where $\Omega_{\infty}=\big\{ x\in\mathbb R^n: p(x)=\infty\big\}$. When $|\Omega_{\infty}|=0$, it is straightforward
\[
F_{p}(f/\eta):=\int\limits_{\mathbb R^n}\left({\frac{|f(x)|}{\eta}}\right)^{p(x)}dx<\infty.
\]
\end{definition}
The variable exponent Lebesgue space $L^{p(\cdot)}(\mathbb R^n)$ then becomes a norm space equipped with a norm as follows
$$\left\|f\right\|_{L^{p(\cdot)}}= \inf \left\{\eta>0: F_p\left(\frac{f}{\eta}\right)\leq 1\right\}.$$
Let us denote by $\mathcal P_{b}(\mathbb R^n)$  the class of exponents $q(\cdot)\in\mathcal P(\mathbb R^n)$ such that
\[
1 < q_{-}\leq q(x) \leq q_{+}<\infty,\,  \text{ for all }\,x\in\mathbb R^n,
\]
where $q_{-}= \text{ess\,inf}_{x\in\mathbb R^n}q(x)$ and $q_{+}= \text{ess\,sup}_{x\in\mathbb R^n}q(x)$.
For $p\in\mathcal P_{b}(\mathbb R^n)$, it is useful to remark that we have the following inequalities which are usually used in the sequel.
\begin{eqnarray}\label{maxminvar}
&&\hskip -8pt [i]\,\, \textit{ \rm If }\, F_{p}(f)\leq C, \textit{ \rm then}\, \big\|f\big\|_{L^{p(\cdot)}} \leq \textit{\rm max}\big\{ {C}^{\frac{1}{q_-}}, {C}^{\frac{1}{q_+}} \big\},\,\textit{\rm for all}\,\, f\in L^{p(\cdot)},
\nonumber
\\
&&\hskip - 8pt [ii]\textit{ \rm If }\, F_{p}(f)\geq C, \textit{ \rm then}\, \big\|f\big\|_{L^{p(\cdot)}}\geq \textit{\rm min}\big\{ {C}^{\frac{1}{q_-}}, {C}^{\frac{1}{q_+}} \big\},\,\textit{\rm for all}\, f\in L^{p(\cdot)}.
\end{eqnarray}
\\
The space ${\mathcal P}_\infty(\mathbb R^n)$ is defined by the set of all measurable functions $q(\cdot)\in\mathcal P(\mathbb R^n)$ and there exists a constant $q_\infty$ such that
\[
q_{\infty} =\lim\limits_{|x|\to\infty}{q(x)}.
\]
\\
For $p(\cdot)\in\mathcal P(\mathbb R^n)$, the weighted variable exponent Lebesgue space $L^{p(\cdot)}_\omega(\mathbb R^n)$ is the set of all complex-valued measurable functions $f$ such that $f\omega$ belongs the $L^{p(\cdot)}(\mathbb R^n)$ space and $f$ has norm
\[
\big\|f\big\|_{L^{p(\cdot)}_\omega}=\big\|f\omega\big\|_{L^{p(\cdot)}}.
\]
Let $\mathbf C_0^{\text{log}}(\mathbb R^n) $ denote the set of  all log-H\"{o}lder continuous functions $\alpha (\cdot)$ satisfying at the origin
\[
\left|\alpha(x)-\alpha(0)\right|\leq \dfrac{C_0^\alpha}{\text{log}\left(e+\frac{1}{|x|}\right)},\,  \text{ for all }\, x\in\mathbb R^n.
\]
Denote by $\mathbf C_\infty^{\text{log}}(\mathbb R^n) $  the set of all log-H\"{o}lder continuous functions $\alpha (\cdot)$ satisfying at infinity
\[
\left|\alpha(x)-\alpha_{\infty}\right|\leq \dfrac{C_\infty^\alpha}{\text{log}(e+|x|)},\,  \text{ for all }\, x\in\mathbb R^n,
\]

Next, we  give the definition of variable exponent weighted Herz spaces ${\mathop{K}\limits^.}^{\alpha(\cdot),p}_{q(\cdot),\omega}$ and variable exponent weighted Morrey-Herz spaces ${M\mathop{K}\limits^.}^{\alpha(\cdot),\lambda}_{p,q(\cdot),\omega}$(see \cite{LZ2014}, \cite{WZ2016} for more details) .
\begin{definition}
Let $0<p<\infty,  q(\cdot)\in \mathcal P_b(\mathbb R^n)$ and $\alpha(\cdot):\mathbb R^n\to \mathbb R$ with $\alpha(\cdot)\in L^{\infty}(\mathbb R^n)$. The variable exponent weighted Herz space ${\mathop{K}\limits^.}^{\alpha(\cdot),p}_{q(\cdot),\omega}$ is defined
by
\[
{\mathop{K}\limits^.}^{\alpha(\cdot),p}_{q(\cdot),\omega}=\left\{f\in L^{q(\cdot)}_\text {loc}(\mathbb R^n\setminus\{0\}):\|f\|_{{\mathop{K}\limits^.}^{\alpha(\cdot),p}_{q(\cdot),\omega}} <\infty \right \},
\]
where $\|f\|_{{\mathop{K}\limits^.}^{\alpha(\cdot),p}_{q(\cdot),\omega}}=\left(\sum\limits_{k=-\infty}^{\infty} \|2^{k\alpha(\cdot)} f\chi_k\|^p_{L^{q(\cdot)}_\omega}\right)^{\frac{1}{p}}.$
\end{definition}
\begin{definition}
Assume that $0\leq \lambda<\infty, 0<p<\infty,  q(\cdot)\in \mathcal P_b(\mathbb R^n)$ and $\alpha(\cdot):\mathbb R^n\to \mathbb R$ with $\alpha(\cdot)\in L^{\infty}(\mathbb R^n)$. The variable exponent weighted Morrey-Herz space ${M\mathop{K}\limits^.}^{\alpha(\cdot),\lambda}_{p,q(\cdot), \omega}$ is defined
by
\[
{M\mathop{K}\limits^.}^{\alpha(\cdot),\lambda}_{p,q(\cdot),\omega}=\left\{f\in L^{q(\cdot)}_\text {loc}(\mathbb R^n\setminus\{0\}):\|f\|_{{M\mathop{K}\limits^.}^{\alpha(\cdot),\lambda}_{p,q(\cdot),\omega}} <\infty \right \},
\]
where $\|f\|_{{M\mathop{K}\limits^.}^{\alpha(\cdot),\lambda}_{p,q(\cdot),\omega}}=\sup\limits_{k_0\in\mathbb Z}2^{-k_0\lambda}\left(\sum\limits_{k=-\infty}^{k_0} \|2^{k\alpha(\cdot)} f\chi_k\|^p_{L^{q(\cdot)}_\omega}\right)^{\frac{1}{p}}.$
\end{definition}
It is easy to see that ${M\mathop{K}\limits^.}^{\alpha(\cdot),0}_{p,q(\cdot)}(\mathbb R^n)$ =${\mathop{K}\limits^.}^{\alpha(\cdot), p}_{q(\cdot)}(\mathbb R^n)$. Consequently, the Herz space with variable exponent is a special case of Morrey-Herz space with variable exponent.
\vskip 5pt
Let us next state the following corollary which is used in the sequel.  The proof is trivial and may be found in \cite{WZ2016}.
\begin{lemma}\label{lemmaVE}
Let $\alpha(\cdot)\in L^{\infty}(\mathbb R^n)$, $q(\cdot)\in \mathcal P_b(\mathbb R^n)$, $p\in(0, \infty)$ and $\lambda\in [0,\infty)$. If $\alpha(\cdot)$ is log-H\"{o}lder continuous
both at the origin and at infinity, then
\[
\big\|f\chi_j\big\|_{L^{q(\cdot)}_\omega} \leq C. 2^{j(\lambda-\alpha(0))}\big\|f\big\|_{{M\mathop{K}\limits^.}^{\alpha(\cdot),\lambda}_{p,q(\cdot),\omega}},\,\,\textit{for all}\,\, j\in\mathbb Z^-,
\]
and
\[
\big\|f\chi_j\big\|_{L^{q(\cdot)}_\omega} \leq C. 2^{j(\lambda-\alpha_\infty)}\big\|f\big\|_{{M\mathop{K}\limits^.}^{\alpha(\cdot),\lambda}_{p,q(\cdot),\omega}},\,\,\textit{for all}\,\, j\in\mathbb N.
\]
\end{lemma}
We recall the definition of two-weight $\lambda$-central Morrey spaces with variable-exponent (see \cite {HausdoffCDD}).
\begin{definition}
For $\lambda \in \mathbb R$ and $p\in\mathcal P_\infty(\mathbb R^n)$, we denote ${\mathop B\limits^.}_{\omega_1,\omega_2}^{p(\cdot),\lambda}$ the class of locally integrable functions $f$ on $\mathbb R^n$ satisfying 
\[
{\big\| f \big\|_{{\mathop B\limits^.}_{\omega_1,\omega_2}^{p(\cdot),\lambda }}} = \mathop {\sup }\limits_{R\, > 0}  {\frac{1}{{{\omega_1}\big(B(0,R)\big)^{\lambda+\frac{1}{p_{\infty}}}}}} \big\|f\big\|_{L^{p(\cdot)}_{\omega_2}(B(0,R))}< \infty,
\]
where $\big\|f\big\|_{L^{p(\cdot)}_{\omega_2}(B(0,R))}=\big\|f\chi_{B(0,R)}\big\|_{L^{p(\cdot)}_{\omega_2}}$
and $\omega_1$, $\omega_2$ are non-negative and local integrable functions.  
\end{definition}
Next, the following theorem is stated as the embedding result on Lebesgue spaces with variable exponent (see \cite {HausdoffCDD}).
\begin{theorem}\label{theoembed}
Let $p(\cdot), q(\cdot)\in\mathcal P(\mathbb R^n)$ and $q(x)\leq p(x)$ almost everywhere $x\in \mathbb R^n$, and
\[
\frac{1}{r(\cdot)}:=\frac{1}{q(\cdot)}-\frac{1}{p(\cdot)}\,\,\textit{\rm and}\,\, \big\|1\big\|_{L^{r(\cdot)}}<\infty.
\]
Then there exists a constant K such that
\[
\big\|f\big\|_{L^{q(\cdot)}_\omega}\leq K \big\|1\big\|_{L^{r(\cdot)}}\big\|f\big\|_{L^{p(\cdot)}_\omega}.
\]
\end{theorem}
Let us recall to define Lipschitz space and central BMO space (see, for example, \cite{JN1961}, \cite{LY1993}, \cite{Stein} for more details).
\begin{definition}
Let $0<\beta\leq 1$. The Lipschitz space $\rm Lip^\beta(\mathbb R^n)$ is defined as the
set of all functions $f:\mathbb R^n\to \mathbb C$ satisfying $\big\|f\big\|_{\rm Lip^\beta(\mathbb R^n)}< \infty$, where
\[
\big\|f\big\|_{\rm Lip^\beta(\mathbb R^n)} := \mathop {\rm sup}\limits_{x, y\in \mathbb R^n, x\neq y}
\dfrac{|f(x)-f(y)|}{|x-y|^\beta}.
\]
\end{definition}
\begin{definition} 
The space $BMO(\mathbb R^n)$ consists of all locally integrable functions $f:\mathbb R^n\to \mathbb C$  satisfying 
\[
\big\|f\big\|_{\rm BMO(\mathbb R^n)} := \mathop {\rm sup}\limits_{Q}\dfrac{1}{|Q|}\int\limits_{Q}{|f(x)-f_{Q}|dx}<\infty,
\]
where the supremum is taken over all cubes $Q \subset \mathbb R^n$ with sides parallel to the coordinate axes.
\end{definition}
\begin{definition}
Let $1\leq q<\infty$ and $\omega$ be a weight function. The central bounded mean oscillation space ${\mathop {CMO}\limits^{.}}^q(\omega)$ is defined as the set of all functions $f\in L^q_{\rm loc}(\mathbb R^n)$ such that
\[
\big\|f\big\|_{{\mathop {CMO}\limits^{.}}^q(\omega)}=\mathop {\rm sup}\limits_{R >0}\Big( \frac{1}{\omega(B(0,R))}\int\limits_{B(0,R)}{|f(x)-f_{\omega, B(0,R)}|^q\omega(x)dx}\Big)^{\frac{1}{q}},
\]
where 
\[
\omega(B(0,R))=\int\limits_{B(0,R)}{\omega(x)dx}\,\,\textit{\rm and}\,\,f_{\omega, B(0,R)}=\frac{1}{\omega(B(0,R))}\int\limits_{B(0,R)}{f(x)\omega(x)dx}.
\]
\end{definition} 
Remark that, Fefferman \cite{Stein} obtain the famous result that the space $BMO(\mathbb R^n)$ is the dual space of Hardy space $H^1(\mathbb R^n)$. When $\omega =1$, we write simply $CMO^q(\mathbb R^n):=CMO^q(\omega)$. The space $CMO(\mathbb R^n)$ can be seen as a local version of $BMO(\mathbb R^n)$ at the origin. Moreover, $BMO(\mathbb R^n)\subsetneqq CMO^q(\mathbb R^n)$, where $1 \leq q < \infty$, and  the John-Nirenberg inequality is not true in $CMO^q(\mathbb R^n)$.
\section{Statement of the results}
Before stating our main results, we introduce some notations which will be used throughout this section. Let $ \gamma_1,...,\gamma_m\in\mathbb R,\lambda_1,...,\lambda_m \geq 0, p_1, ..., p_m, p\in (0,\infty)$, $0<\beta_1,...,\beta_m \leq 1$, $q_i\in\mathcal P_b(\mathbb R^n), r_i\in\mathcal P_\infty(\mathbb R^n)$ for $i=1,...,m$ and $\alpha_1,...,\alpha_m\in L^\infty(\mathbb R^n)\cap \mathbf C_0^{\text{log}}(\mathbb R^n)\cap \mathbf C_{\infty}^{\text{log}}(\mathbb R^n)$. The functions $\alpha^*(\cdot),q(\cdot),\gamma(\cdot)$ and numbers $\beta,\lambda$ are defined as follows
$$\beta_1+\cdots+\beta_m=\beta,$$
$$  {\lambda _1} + {\lambda _2} + \cdots + {\lambda _m} = \lambda. $$
$$\gamma_1+\cdots+\gamma_m+\frac{\gamma_1}{r_1(\cdot)}+\cdots+\frac{\gamma_m}{r_m(\cdot)} = \gamma(\cdot),$$
  $$\frac{1}{{{q_1(\cdot)}}}+ \cdots  + \frac{1}{{{q_m(\cdot)}}} +\frac{1}{r_1(\cdot)}+\cdots+\frac{1}{r_m(\cdot)}= \frac{1}{q(\cdot)}, $$
$$ {\alpha _1(\cdot)} + \cdots + {\alpha _m(\cdot)}-\beta_1-\cdots-\beta_m-\frac{\gamma_1+n}{r_1(\cdot)}-\cdots-\frac{\gamma_m+n}{r_m(\cdot)} = \alpha^*(\cdot).$$
For a matrix $A=(a_{ij})_{n\times n}$, we define the norm of $A$ as follow
\begin{equation}\label{normB}
\left\|A\right\|=\left(\sum\limits_{i,j=1}^{n}{{|a_{ij}|}^2}\right)^{1/2}.
\end{equation}
As above we conclude that $\left|Ax\right|\leq \left\|A\right\|\left|x\right|$ for any vector $x\in\mathbb R^n$. In particular, if $A$ is invertible,  then we have
\begin{equation}\label{detA}
\left\|A\right\|^{-n}\leq \left|\rm det(A^{-1})\right|\leq \left\|A^{-1}\right\|^{n}.
\end{equation}
Now, we are ready to state the main results in this paper.
\begin{theorem}\label{TheoremMorreyHerzLip}
Let $\zeta >0$, $\omega_1(x)=|x|^{\gamma_1}, ... ,  \omega_m(x)=|x|^{\gamma_m}, \omega(x)=|x|^{\gamma(x)}$, $q(\cdot)\in \mathcal P_b(\mathbb R^n)$, $\alpha^*\in L^\infty(\mathbb R^n)\cap \mathbf C_0^{\text{log}}(\mathbb R^n)\cap \mathbf C_{\infty}^{\text{log}}(\mathbb R^n)$, $b_i\in {\rm Lip}^{\beta_i}$, $\lambda_1,...,\lambda_m>0$ and the following conditions are true:
\begin{equation}\label{DKnhung}
q_i(A_i^{-1}(t)\cdot)\leq \zeta.q_i(\cdot)\,\textit{\rm and}\,\big\|1\big\|_{L^{\vartheta_{i}(t,\cdot)}}<\infty, \, \text{a.e. }\,t\in \textit{\rm supp}(\Phi),\,\textit{\rm for all}\,\,i=1,...,m,
\end{equation}
\begin{equation}\label{DKalpha}
\alpha_i(0)-\alpha_{i\infty}\geq 0,\,\,\textit{\rm for all}\,\,i=1,...,m,
\end{equation}
either
$$
\gamma_1,...,\gamma_m>-n, r_1(0)= r_{1+}, r_{1\infty}=r_{1-},..., r_m(0)= r_{m+},  r_{m\infty}=r_{m-}\, $$
or
$$\gamma_1,...,\gamma_m<-n\,, r_1(0)= r_{1-}, r_{1\infty}=r_{1+}, ..., r_m(0)= r_{m-},  r_{m\infty}=r_{m+},$$
or
\begin{eqnarray}\label{DKgammari}
\gamma_1=\cdots=\gamma_m=-n.
\end{eqnarray}
Then, if
\begin{equation}\label{Var-MorreyHerzDK4}
\mathcal C_1=\int\limits_{\mathbb R^n}{\frac{\Phi(t)}{|t|^n}\prod\limits_{i=1}^{m} c_{A_i,q_i,\gamma_i}(t)}\big\|I_n-A_i(t)\big\|^{\beta_i}\big\|1\big\|_{L^{\vartheta_{i}(t,\cdot)}}\phi_{A_i,\lambda_i}(t)dt<\infty,
\end{equation}
where 
\begin{eqnarray}
\phi_{A_i,\lambda_i}(t)&=&{\max\Big\{\big\|A_i(t)\big\|^{\lambda_i-\alpha_i(0)}, \big\|A_i(t)\big\|^{\lambda_i-\alpha_{i\infty}} \Big\}}\times\nonumber
\\
&&\,\,\,\,\,\,\,\,\,\,\,\,\,\times{\max\Big\{\sum\limits_{r=\Theta_n^*-1}^{0} 2^{r(\lambda_i-\alpha_i(0))},\sum\limits_{r=\Theta_n^*-1}^{0}2^{r(\lambda_i-\alpha_{i\infty})}\Big\}},
\end{eqnarray}
with $\Theta_n^*=\Theta_n^*(t)$ is the greatest integer number satisfying 
\[
\mathop{\rm max}\limits_{i=1,m}\big\{\|A_i(t)\|.\|A_i^{-1}(t)\|\big\}<2^{-\Theta_n^*},\,\, \text{\rm for a.e. }\, t\in\mathbb R^n,
\]
\[
{c_{A_i,q_i,\gamma_i}(t)} = {{{\rm max }}\big\{ {{{\big\| {{A_i}(t)} \big\|}^{ - {\gamma _i}}},{{\big\| {A_i^{ - 1}(t)} \big\|}^{{\gamma _i}}}} \big\}} {\rm max}\big\{{\left| {\det A_i^{ - 1}(t)} \right|}^{\frac{1}{q_{i+}}},{\left| {\det A_i^{ - 1}(t)} \right|}^{\frac{1}{q_{i-}}}\big\},
\]
\[
\dfrac{1}{\vartheta_i(t,\cdot)}=\dfrac{1}{q_i(A_i^{-1}(t)\cdot)}-\dfrac{1}{\zeta q_i(\cdot)},\,\textit{\rm for all} \,i=1,...,m,
\]
we have  $H_{\Phi,\vec{A}}^{\vec b}$ is a bounded operator from ${M{\mathop K\limits^.}}^{\alpha_1(\cdot),\lambda_1}_{p_1,\zeta q_1(\cdot),\omega_1}\times\cdots\times {M{\mathop K\limits^.}}^{\alpha_m(\cdot),\lambda_m}_{p_m,\zeta q_m(\cdot),\omega_m}$ to ${M{\mathop K\limits^.}}^{\alpha^*(\cdot),\lambda}_{p,\zeta q(\cdot),\omega}$.
\end{theorem}
\begin{theorem}\label{TheoremHerzLip}
Suppose that we have the given supposition of Theorem $\ref{TheoremMorreyHerzLip}$. Let $1\leq p, p_1,...,p_m <\infty$, $\lambda_i=0$ and  $\alpha_i(0)=\alpha_{i\infty}$, for all $i=1,...,m$. At the same time, let
\begin{equation}\label{DKpip}
\frac{1}{p_1}+\cdots+\frac{1}{p_m}=\frac{1}{p},
\end{equation}
\begin{equation}\label{DKMorreyHerz} 
\mathcal C_2=\int\limits_{\mathbb R^n}{(2 - {\Theta_n^*})^{m - \frac{1}{p}}\frac{\Phi(t)}{|t|^n}\prod\limits_{i=1}^{m}}c_{A_i,q_i,\gamma_i}(t)\big\|I_n-A_i(t)\big\|^{\beta_i}\big\|1\big\|_{L^{\vartheta_{i}(t,\cdot)}}\phi_{A_i,0}(t)dt<\infty,
\end{equation}
Then, $H_{\Phi,\vec{A}}^{\vec b}$ is a bounded operator from ${\mathop K\limits^.}_{\zeta q_1(\cdot),\omega_1}^{{\alpha _1(\cdot)}, p_1}\times \cdots\times {\mathop K\limits^.}_{\zeta q_m(\cdot),\omega_m}^{{\alpha _m(\cdot)}, p_m}$ to ${\mathop K\limits^.}_{q(\cdot),\omega}^{{\alpha^*(\cdot)}, p}.$
\end{theorem}
By using the ideas in the proof of Theorem \ref{TheoremMorreyHerzLip}, we give the analogous result for the Lebesgue spaces with variable exponent as follows.
\begin{theorem}\label{TheoremLebLip}
Let $\zeta >0$, $\gamma_1,...,\gamma_m<0$, $\omega_1(x)=|x|^{\gamma_1},  ... ,  \omega_m(x)=|x|^{\gamma_m}, \omega(x)=|x|^{\gamma(x)}$, $q(\cdot)\in \mathcal P_b(\mathbb R^n)$, $b_1,...,b_m\in {\rm Lip}^{\beta_i}$, and let the hypothesis (\ref{DKnhung}) in Theorem \ref{TheoremMorreyHerzLip} hold. Thus, if the following conditions are true:
\begin{equation}\label{DK|x|^betai}
\big\||\cdot|^{\beta_i+\frac{\gamma_i}{r_i(\cdot)}}\big\|_{L^{r_i(\cdot)}_{\omega_i}}<\infty,\,\textit{\rm for all}\,\,i=1,...,m,
\end{equation}
\begin{equation}\label{Var-MorreyHerzDK4}
\mathcal C_3=\int\limits_{\mathbb R^n}{\frac{\Phi(t)}{|t|^n}\prod\limits_{i=1}^{m}}c_{A_i,q_i,\gamma_i}(t).\big\|I_n-A_i(t)\big\|^{\beta_i}\big\|1\big\|_{L^{\vartheta_{i}(t,\cdot)}}dt<\infty,
\end{equation}
then we have
 \[
\big\|H_{\Phi,\vec A}^{\vec b}(\vec f) \big\|_{L^{q(\cdot)}_{\omega}}\lesssim \mathcal C_3.\mathcal B_{\rm Lip}.\prod\limits_{i = 1}^m {{{\left\| {{f_i}} \right\|}_{L^{\zeta q_i(\cdot)}_{\omega_i}}}}.
 \]
\end{theorem}
Next, we consider that all of $r_1(\cdot), ..., r_m(\cdot)$ are constant and the following conditions hold:
\begin{enumerate}
\item[$(\rm H_1)$] $\alpha_1(\cdot)+\cdots+\alpha_m(\cdot)-\frac{\gamma_1+n}{r_1}-\cdots-\frac{\gamma_m+n}{r_m}=\alpha^{**}(\cdot)$,
\\
\item[$(\rm H_2)$] $A_i(t)= s_i(t).a_i(t)$ for all $i=1, ..., m$, where $s_i: \rm supp(\Phi)\to \mathbb R$ is a measurable function such that $s_i(t)\neq 0$ for a.e  $t\in \rm supp(\Phi)$ and $a_i(t)$ is an $n\times n$ rotation matrix for a.e  $t\in \rm supp(\Phi)$.
\end{enumerate}
Then, we also obtain the following some interesting results.
\begin{theorem}\label{TheoremMorreyHerzCMO}
Let $\zeta >0$, $\lambda_1, ..., \lambda_m>0$, $\gamma_1, ..., \gamma_m>-n$, $\omega_1(x)=|x|^{\gamma_1}, ... $, $\omega_m(x)=|x|^{\gamma_m}, \omega(x)=|x|^{\gamma(x)}$, $q(\cdot)\in \mathcal P_b(\mathbb R^n)$, $b_1\in{\mathop{CMO}\limits^{.}}^{r_1}(\omega_1)$,$ ...$, $b_m\in {\mathop {CMO}\limits^{.}}^{r_m}(\omega_m)$, the hypothesis (\ref{DKnhung}) and (\ref{DKalpha}) in Theorem \ref{TheoremMorreyHerzLip} hold. Then, if 
\begin{equation}\label{DKMorreyHerzCMO}
\mathcal C_4 =\int\limits_{\mathbb R^n}{\frac{\Phi(t)}{|t|^n}\prod\limits_{i=1}^{m}}c_{A_i,q_i,\gamma_i}(t).\big\|1\big\|_{L^{\vartheta_{i}(t,\cdot)}}\phi_{A_i,\lambda_i}(t)\big(1+\psi_{A_i,\gamma_i}^{\frac{1}{r_i}}|s_i(t)|^{\frac{\gamma_i+n}{r_i}}+\varphi_{A_i}(t)\big)dt<\infty,
\end{equation}
where
\[
\psi_{A_i,\gamma_i}(t)=|{\rm det} A_i^{-1}(t)|{\rm max}\big\{\big\|A_i^{-1}(t)\big\|^{\gamma_i},\big\|A_i(t)\big\|^{-\gamma_i}\big\},
\]
\[
\varphi_{A_i}(t)={\rm max}\big\{{\rm log}(4|s_i(t)|),{\rm log}\frac{2}{|s_i(t)|}\big\},
\]
we have $H_{\Phi,\vec{A}}^{\vec b}$ is a bounded operator from ${M{\mathop K\limits^.}}^{\alpha_1(\cdot),\lambda_1}_{p_1,\zeta q_1(\cdot),\omega_1}\times\cdots\times {M{\mathop K\limits^.}}^{\alpha_m(\cdot),\lambda_m}_{p_m,\zeta q_m(\cdot),\omega_m}$ to ${M{\mathop K\limits^.}}^{\alpha^{**}(\cdot),\lambda}_{p,\zeta q(\cdot),\omega}$.
\end{theorem}
\begin{theorem}\label{TheoremHerCMO}
Let  $1\leq p, p_1,...,p_m <\infty$, $\lambda_i=0$, $\alpha_i(0)=\alpha_{i\infty}$, for all $i=1,...,m$. Also, both the assumptions of Theorem \ref{TheoremMorreyHerzCMO} and the hypothesis (\ref{DKpip}) in Theorem \ref{TheoremHerzLip}. In addition, the following condition holds:
\begin{eqnarray}\label{DKHerzCMO}
\mathcal C_5 &=&\int\limits_{\mathbb R^n}{(2 - {\Theta_n^*})^{m - \frac{1}{p}}\frac{\Phi(t)}{|t|^n}\prod\limits_{i=1}^{m}}c_{A_i,q_i,\gamma_i}(t)\big\|1\big\|_{L^{\vartheta_{i}(t,\cdot)}}\phi_{A_i,0}(t)\times
\nonumber
\\
&&\,\,\,\,\,\,\,\,\,\,\,\,\,\,\,\,\,\,\,\times\big(1+\psi_{A_i,\gamma_i}^{\frac{1}{r_i}}|s_i(t)|^{\frac{\gamma_i+n}{r_i}}+\varphi_{A_i}(t)\big)dt<\infty.
\end{eqnarray}
Then, we have
\[
\big\|H_{\Phi,\vec{A}}^{\vec b}(\vec f)\big\|_{{\mathop K\limits^.}_{q(\cdot),\omega}^{{\alpha^{**}(\cdot)}, p}}\lesssim \mathcal C_5\prod\limits_{i=1}^m \big\|f_i\big\|_{{\mathop K\limits^.}_{\zeta q_i(\cdot),\omega_i}^{{\alpha _i(\cdot)}, p_i}}.
\]
\end{theorem}
Let us now assume that $q(\cdot)$ and $q_i(\cdot)\in\mathcal P_\infty(\mathbb R^n)$, $\lambda,\alpha,\gamma,\beta, r_i, {\lambda _i},\alpha_i, {\gamma_i}, \beta_i$ are real numbers such that  $r_i\in (0,\infty)$, $\lambda_i\in \big(\frac{-1}{q_{i\infty}},0\big)$, $\alpha_i, \gamma_i \in (-n,\infty)$, $\beta_i\in (0,1]$, $i=1,2,...,m$ and denote
  $$\beta_1+\cdots+\beta_m=\beta,$$
$$ \alpha_1+\cdots+\alpha_m+\frac{\alpha_1}{r_1}+\cdots+\frac{\alpha_m}{r_m}=\alpha,$$
$$   \frac{1}{{{q_1}(\cdot)}}+ \cdots + \frac{1}{{{q_m}(\cdot)}}+\frac{1}{r_1} +\cdots + \frac{1}{r_m} = \frac{1}{q(\cdot)}. $$
We are also interested in the commutators of multilinear Hausdorff operators on the product of weighted $\lambda$-central Morrey spaces with variable exponent. More precisely, we have the following useful result.
\begin{theorem}\label{TheoremMorreyLip}
Let $\omega_i(x)=|x|^{\gamma_i}, v_i(x)=|x|^{\alpha_i}, b_i\in {\rm Lip}^{\beta_i}$ for all $i= 1,...,m$, $\omega(x)= |x|^\gamma,v(x)=|x|^{\alpha}$ and  the following conditions are true: 
\begin{equation}\label{DKnhung1}
q_i(A_i^{-1}(t)\cdot)\leq q_i(\cdot)\,\textit{\rm and}\,\big\|1\big\|_{L^{\vartheta_{1i}(t,\cdot)}}<\infty, \, \text{a.e. }\,t\in \textit{\rm supp}(\Phi),\,\textit{\rm for all}\,\,i=1,...,m,
\end{equation}
\begin{equation}\label{DKlambda}
\beta+\alpha-\frac{\gamma}{q_{\infty}}+ \sum\limits_{i=1}^m (\gamma_i+n)\lambda_i-\alpha_i+
\frac{\gamma_i}{q_{i\infty}}=(\gamma+n)\lambda,
\end{equation}
\begin{eqnarray}\label{DK2MVar}
\mathcal C_{6} &=& \int\limits_{\mathbb R^n} { \frac{{ {\Phi (t)}}}{{{{\left| t \right|}^n}}}\prod\limits_{i = 1}^m {{{\left\| {A_i(t)} \right\|}^{(n + {\gamma _i})\left(\frac{1}{q_{i\infty}}+{\lambda _i}\right)}}}}c_{A_i,q_i,\alpha_i}(t)\times
\nonumber
\\
&&\,\,\,\,\,\,\,\,\,\,\,\,\,\,\,\,\,\,\,\,\,\,\,\,\,\,\,\,\,\,\;\;\;\;\;\;\;\;\times\big\|1\big\|_{L^{\vartheta_{1i}(t,\cdot)}}\big\|I_n-A_i(t)\big\|^{\beta_i}dt <+ \infty,
\end{eqnarray}
where 
\[
\dfrac{1}{\vartheta_{1i}(t,\cdot)}=\dfrac{1}{q_i(A_i^{-1}(t)\cdot)}-\dfrac{1}{q_i(\cdot)},\,\textit{\rm for all} \,i=1,...,m.
\]
Then, we have $H_{\Phi ,\vec{A}}^{\vec b}$ is bounded from ${\mathop B\limits^.}_{\omega_1,v_1}^{q_1(\cdot),\lambda _1}\times  \cdots\times {\mathop B\limits^.}_{{\omega_m,v_m}}^{{q_m(\cdot)},{\lambda _m}}$ to ${\mathop B\limits^.}_{\omega,v}^{q(\cdot),\lambda}$. 
\end{theorem}
\begin{theorem}\label{TheoremMorreyCMO}
Given $\omega_i(x)=|x|^{\gamma_i},v_i(x)=|x|^{\alpha_i}, b_i\in {\mathop{CMO}\limits^{.}}^{r_i}(\omega_i)$ for all $i= 1,...,m$, $\omega(x)= |x|^\gamma, v(x)=|x|^{\alpha}$, the hypothesis (\ref{DKnhung1}) in Theorem \ref{TheoremMorreyLip} and the condition ($H_2$) hold. In addition, the following statements are true:
\begin{equation}\label{DKlambda1}
\alpha-\frac{\gamma}{q_{\infty}}+ \sum\limits_{i=1}^m (\gamma_i+n)\lambda_i-\alpha_i+
\frac{\gamma_i}{q_{i\infty}}=(\gamma+n)\lambda,
\end{equation}
\begin{eqnarray}\label{DK2MVar}
\mathcal C_{7}&=& \int\limits_{\mathbb R^n} { \frac{{ {\Phi (t)}}}{{{{\left| t \right|}^n}}}\prod\limits_{i = 1}^m {{{\left\| {A_i(t)} \right\|}^{(n + {\gamma _i})\left(\frac{1}{q_{i\infty}}+{\lambda _i}\right)}}}}c_{A_i,q_i,\alpha_i}(t).\big\|1\big\|_{L^{\vartheta_{1i}(t,\cdot)}}\times\nonumber
\\
&&\,\,\,\,\,\,\,\,\,\,\,\,\,\,\,\,\,\,\,\,\,\,\,\,\,\,\,\,\,\,\,\,\,\times \big(1+\psi_{A_i,\alpha_i}^{\frac{1}{r_i}}|s_i(t)|^{\frac{\alpha_i+n}{r_i}}+\varphi_{A_i}(t)\big)dt <+ \infty.
\end{eqnarray}
Then, we conclude that $H_{\Phi ,\vec{A}}^{\vec b}$ is bounded from ${\mathop B\limits^.}_{\omega_1,v_1}^{q_1(\cdot),\lambda _1}\times  \cdots\times {\mathop B\limits^.}_{{\omega_m},v_m}^{{q_m(\cdot)},{\lambda _m}}$ to ${\mathop B\limits^.}_{\omega,v}^{q(\cdot),\lambda}$. 
\end{theorem}
\section{Proofs of the theorems}
Fristly, for simplicity of notation, we denote 
\[
\mathcal B_{\rm Lip}=\prod\limits_{i=1}^m\big\|b_i\big\|_{{\rm Lip}^{\beta_i}},\,\mathcal B_{{\rm CMO},\vec{\omega}}=\prod\limits_{i=1}^m\big\|b_i\big\|_{{\mathop {CMO}\limits^.}^{r_i}(\omega_i)}\,\textit{\rm and}\, \mathcal F=\prod\limits_{i=1}^m\big\|f_i\big\|_{{M{\mathop K\limits^.}}^{\alpha_i(\cdot),\lambda_i}_{p_i,\zeta q_i(\cdot),\omega_i}}.
\]
\subsection{Proofs of Theorem \ref {TheoremMorreyHerzLip} and Theorem \ref{TheoremHerzLip}}

By using the versions of the Minkowski inequality for variable Lebesgue spaces from Corollary 2.38 in \cite{CUF2013}, we have
\begin{equation}\label{Hvec}
\big\|H_{\Phi,\vec{A}}^{\vec b}(\vec{f})\chi_k\big\|_{L^{q(\cdot)}_\omega}\lesssim\int\limits_{\mathbb R^n}{\frac{\Phi(t)}{|t|^n}\big\|\prod\limits_{i=1}^{m}f(A_i(t)\cdot)\big(b_i(\cdot)- b_i(A_i(t)\cdot)\big)\chi_k\big\|_{L^{q(\cdot)}_\omega}}dt.
\end{equation}
On the other hand, since $b_i\in {\rm Lip}^{\beta_i}$, we get
\[
|f(A_i(t)x)\big(b_i(x)- b_i(A_i(t)x)\big)\chi_k(x)|\leq |f(A_i(t)x)|.\big\|b_i\big\|_{{\rm Lip}^{\beta_i}}\big\|I_n-A_i(t)\big\|^{\beta_i}2^{\beta_ik}\chi_k(x).
\]
Thus, by applying the H\"{o}lder inequality for variable Lebesgue spaces (see also Corollary 2.30  in \cite{CUF2013}), we find
\begin{eqnarray}\label{fLip}
&&\big\|\prod\limits_{i=1}^{m}f(A_i(t)\cdot)\big(b_i(\cdot)- b_i(A_i(t)\cdot)\big)\chi_k\big\|_{L^{q(\cdot)}_\omega}
\\
&&\,\,\,\,\,\,\,\,\,\leq 2^{k\beta}\mathcal B_{\rm Lip}\prod\limits_{i=1}^m\big\|I_n-A_i(t)\big\|^{\beta_i}.\prod_{i=1}^{m}\big\|f_i(A_i(t)\cdot)\chi_k\big\|_{L^{q_i(\cdot)}_{\omega_i}}\big\||\cdot|^{\frac{\gamma_i}{r_i(\cdot)}}\chi_k\big\|_{L^{r_i(\cdot)}}.
\nonumber
\end{eqnarray}
We observe that 
\[F_{r_i}(|\cdot|^{\frac{\gamma_i}{r_i(\cdot)}}\chi_k)=\int\limits_{C_k}{|x|^{\gamma_i}dx}=\int\limits_{2^{k-1}}^{2^{k}}{\int\limits_{S^{n-1}}r^{\gamma_i+n-1}d\sigma(x')dr}\lesssim 2^{k(\gamma_i+n)}.
\]
\\
Case 1: $k< 0$. Denote by
\[
\sigma_i = \left\{ \begin{array}{l}
\frac{1}{r_{i+}},\,\textit{\rm if}\, (\gamma_i+n)>0,
\\
\\
\frac{1}{r_{i-}},\,\textit{\rm otherwise}.
\end{array} \right.\]
Case 2: $k\geq 0$. Denote by
\[
\sigma_i = \left\{ \begin{array}{l}
\frac{1}{r_{i-}},\,\textit{\rm if}\, (\gamma_i+n)>0,
\\
\\
\frac{1}{r_{i+}},\,\textit{\rm otherwise}.
\end{array} \right.\]
From this, by (\ref{maxminvar}), we have
\begin{equation}\label{chik}
\big\||\cdot|^{\frac{\gamma_i}{r_i(\cdot)}}\chi_k\big\|_{L^{r_i(\cdot)}}\lesssim 2^{k(\gamma_i+n)\sigma_i},
\end{equation}
Therefore, from (\ref{Hvec})-(\ref{chik}), we see that
\begin{eqnarray}\label{Hvec1}
\big\|H_{\Phi,\vec{A}}^{\vec b}(\vec{f})\chi_k\big\|_{L^{q(\cdot)}_{\omega}}&\lesssim& 2^{k(\beta+\sum\limits_{i=1}^{m} (\gamma_i+n)\sigma_i)}\mathcal B_{\rm Lip}\times
\\
&&\times \int\limits_{\mathbb R^n}{\frac{\Phi(t)}{|t|^n}}\prod_{i=1}^{m}\big\|I_n-A_i(t)\big\|^{\beta_i}.\big\|f_i(A_i(t)\cdot)\chi_k\big\|_{L^{q_i(\cdot)}_{\omega_i}}dt.
\nonumber
\end{eqnarray}
Let us now fix $i\in \big\{1,2,...,m\big\}$. Since  $\|A_i(t)\|\neq 0$, there exists an integer number $\ell_i=\ell_i(t)$ such that $2^{\ell_i-1}<\|A_i(t)\|\leq 2^{\ell_i}$. 
By writing $\rho_{\vec A}^* (t)$ as
\[
 \rho_{\vec A}^* (t) =\mathop{\rm max}\limits_{i=1,...,m}\big\{\big\|A_i(t)\big\|.\big\|A_i^{-1}(t)\big\|\big\}.
 \]
Hence, by letting $y=A_i(t).z$ with $z\in C_k$, we arrive at
\[
| y | \geq {\left\| {A_i^{ - 1}(t)} \right\|^{ - 1}}\left| z \right|\geq \frac{{{2^{{\ell_i} + k - 2}}}}{{{\rho _{\vec{A}}^*}}} > {2^{k + {\ell_i} - 2 + {\Theta_n^*}}},
\]
and
\[
| y |\leq \left\| {{A_i}(t)} \right\|.\left| z \right| \leq {2^{{\ell_i} + k}}.
\]
These estimations can be used to imply that
\begin{equation}\label{AiCk}
{A_i}(t).{C_k} \subset \left\{ {z \in {\mathbb R^n}:{2^{k + {\ell_i} - 2 + {\Theta _n^*}}} < \left| z \right| \leq {2^{k + {\ell_i}}}} \right\}.
\end{equation}
Now, we will prove the following inequality
\begin{equation}\label{fAchikLq}
\big\|f_i(A_i(t)\cdot)\chi_k\big\|_{L^{q_i(\cdot)}_{\omega_i}}\lesssim c_{A_i,q_i,\gamma_i}(t).\big\|1 \big\|_{L^{r_i(t,\cdot)}}.\sum\limits_{r=\Theta_n^*-1}^{0} \big\|f_i\chi_{k+\ell_i+r}\big\|_{L^{\zeta q_i(\cdot)}_{\omega_i}}.
\end{equation}
Indeed, for $\eta>0$, by (\ref{AiCk}), we get
\begin{eqnarray}
&&\int\limits_{\mathbb R^n}{\left(\dfrac{\big|f_i(A_i(t)x)\chi_k(x)\omega_i(x)\big|}{\eta}\right)^{q_i(x)}dx}\nonumber
\\
&&\leq\int\limits_{A_i(t)C_k}{\left(\dfrac{\big|f_i(z)\big|{\rm max}\big\{\big\|A_i^{-1}(t)\big\|^{\gamma_i},\big\|A_i(t)\big\|^{-\gamma_i}\big\}\omega_i(z)}{\eta}\right)^{q_i(A_i^{-1}(t)z)}\big|
\textit{\rm det} A_i^{-1}(t)\big|dz}\nonumber
\\
&&\leq\int\limits_{\mathbb R^n}\left(\frac{c_{A_i,q_i,\gamma_i}(t)\big|\sum\limits_{r=\Theta_n^*-1}^0 f_i(z)\chi_{k+\ell_i+r}(z)\big|\omega_i(z)}{\eta}\right)^{q_i(A_i^{-1}(t).z)}dz.
\nonumber
\end{eqnarray}
From this, by the definition of Lebesgue space with variable exponent, we find
\[
\big\|f_i(A_i(t)\cdot)\chi_k\big\|_{L^{q_i(\cdot)}_{\omega_i}}\leq c_{A_i,q_i,\gamma_i}(t).\sum\limits_{r=\Theta_n^*-1}^{0} \big\|f_i\chi_{k+\ell_i+r}\big\|_{L^{q_i(A_i^{-1}(t)\cdot)}_{\omega_i}}.
\]
In view of ($\ref{DKnhung}$) and Theorem \ref{theoembed}, we deduce
\[
\big\|f\big\|_{L^{q_i(A_i^{-1}(t)\cdot)}_{\omega_i}}\lesssim \big\|1\big\|_{L^{\vartheta_i(t,\cdot)}}.\big\|f\big\|_{L^{\zeta q_i(\cdot)}_{\omega_i}}.
\]
This completes the proof of the inequalities (\ref{fAchikLq}). 
Now, combining (\ref{Hvec1}) and (\ref{fAchikLq}), it is easy to see that
\begin{eqnarray}\label{Hvec2}
\big\|H_{\Phi,\vec A}^{\vec b}(\vec f)\chi_k\big\|_{L^{q(\cdot)}_{\omega}}&\lesssim& 2^{k(\beta+\sum\limits_{i=1}^{m} (\gamma_i+n)\sigma_i)}\mathcal B_{\rm Lip}\Big(\int\limits_{\mathbb R^n}{\frac{\Phi(t)}{|t|^n}\prod\limits_{i=1}^{m}c_{A_i,q_i,\gamma_i}(t)}\big\|1\big\|_{L^{\vartheta_i(t,\cdot)}}\times
\nonumber
\\
&\times&\big\|I_n-A_i(t)\big\|^{\beta_i}\prod\limits_{i=1}^{m}\sum\limits_{r=\Theta_n^*-1}^{0} \big\|f_i\chi_{k+\ell_i+r}\big\|_{L^{\zeta q_i(\cdot)}_{\omega_i}}dt\Big).
\end{eqnarray}
Thus, by  applying Lemma \ref{lemmaVE} in Section 2, we have
\begin{equation}\label{Hvec3}
\big\|H_{\Phi,\vec A}^{\vec b}(\vec f)\chi_k\big\|_{L^{q(\cdot)}_{\omega}}\lesssim\mathcal B_{\rm Lip}.\mathcal F \Big(\int\limits_{\mathbb R^n}{\frac{\Phi(t)}{|t|^n}\mathcal U(t)\prod\limits_{i=1}^{m}c_{A_i,q_i,\gamma_i}(t)}\big\|1\big\|_{L^{\vartheta_i(t,\cdot)}}\big\|I_n-A_i(t)\big\|^{\beta_i}dt\Big).
\end{equation}
Here
\begin{eqnarray}
\mathcal U(t)&=&2^{k(\beta+\sum\limits_{i=1}^m (\gamma_i+n)\sigma_i)}\prod\limits_{i=1}^{m}\Big(2^{(k+\ell_i)(\lambda_i-\alpha_i(0))}\sum\limits_{r=\Theta_n^*-1}^{0} 2^{r(\lambda_i-\alpha_i(0))}\nonumber
\\
&&\,\,\,\,\,\,\,\,\,\,\,\,\,\,\,\,\,\,\,\,\,\,\,\,\,\,\,\,\,\,\,\,\,\,\,\,\,\,\,\,\,\,\,\,\,\,\,\,\,\,\,\,\,\,+  2^{(k+\ell_i)(\lambda_i-\alpha_{i\infty})}\sum\limits_{r=\Theta_n^*-1}^{0}2^{r(\lambda_i-\alpha_{i\infty})}\Big).
\nonumber
\end{eqnarray}
Since $2^{\ell_i-1}<\big\|A_i(t)\big\|\leq 2^{\ell_i}$, for all $i=1,...,m$, it implies that
\[
2^{\ell_i(\lambda_i-\alpha_i(0))}+2^{\ell_i(\lambda_i-\alpha_{i\infty})}\lesssim \max\big\{\big\|A_i(t)\big\|^{\lambda_i-\alpha_i(0)}, \big\|A_i(t)\big\|^{\lambda_i-\alpha_{i\infty}} \big\}.
\]
From this, we can estimate $\mathcal U$ as follows
\begin{eqnarray}
\mathcal U(t)&\lesssim &2^{k(\beta+\sum\limits_{i=1}^m (\gamma_i+n)\sigma_i)}\prod\limits_{i=1}^{m} {\max\Big\{\big\|A_i(t)\big\|^{\lambda_i-\alpha_i(0)}, \big\|A_i(t)\big\|^{\lambda_i-\alpha_{i\infty}} \Big\}}\times\nonumber
\\
&&\times\Big\{2^{k(\lambda_i-\alpha_i(0))}\sum\limits_{r=\Theta_n^*-1}^{0} 2^{r(\lambda_i-\alpha_i(0))}+ 2^{k(\lambda_i-\alpha_{i\infty})}\sum\limits_{r=\Theta_n^*-1}^{0}2^{r(\lambda_i-\alpha_{i\infty})}\Big\}\nonumber
\\
&\lesssim& 2^{k(\beta+\sum\limits_{i=1}^m (\gamma_i+n)\sigma_i)}\prod\limits_{i=1}^{m} {\max\Big\{\big\|A_i(t)\big\|^{\lambda_i-\alpha_i(0)}, \big\|A_i(t)\big\|^{\lambda_i-\alpha_{i\infty}} \Big\}}\times\nonumber
\\
&&\times{\max\Big\{\sum\limits_{r=\Theta_n^*-1}^{0} 2^{r(\lambda_i-\alpha_i(0))},\sum\limits_{r=\Theta_n^*-1}^{0} 2^{r(\lambda_i-\alpha_{i\infty})}\Big\}}\Big\{2^{k(\lambda_i-\alpha_i(0))}+ 2^{k(\lambda_i-\alpha_{i\infty})}\Big\}.
\nonumber
\end{eqnarray}
This implies that
\[
\mathcal U(t)\lesssim 2^{k(\beta+\sum\limits_{i=1}^m (\gamma_i+n)\sigma_i)} \prod\limits_{i=1}^m\Big\{2^{k(\lambda_i-\alpha_i(0))}+ 2^{k(\lambda_i-\alpha_{i\infty})}\Big\}\phi_{A_i,\lambda}(t).
\]
Thus, by (\ref{Hvec3}), it is not difficult to show that
\begin{equation}\label{Hvec4}
\big\|H_{\Phi,\vec A}^{\vec b}(\vec f)\chi_k\big\|_{L^{q(\cdot)}_{\omega}}\lesssim\mathcal C_1.\mathcal B_{\rm Lip}.\mathcal F. 2^{k(\beta+\sum\limits_{i=1}^m (\gamma_i+n)\sigma_i)}\prod_{i=1}^{m}\big(2^{k(\lambda_i -\alpha_i(0))}+2^{k(\lambda_i -\alpha_{i\infty})}\big).
\end{equation}
Next, using Proposition 2.5 in \cite{LZ2014}, we have
\begin{equation}\label{Hvec5}
\big\| H_{\Phi,\vec A}^{\vec b}(\vec f)\big\|_{M{\mathop K\limits^.}^{\alpha^*(\cdot), \lambda}_{p,q(\cdot),\omega}} \lesssim \max\big\{
\sup\limits_{k_0<0,k_0\in\mathbb Z}E_1, \sup\limits_{k_0\geq 0,k_0\in\mathbb Z}(E_2+E_3)
\big\},
\end{equation}
where 
\begin{eqnarray}
&&E_1=2^{-k_0\lambda}\Big(\sum\limits_{k=-\infty}^{k_0}2^{k\alpha^*(0)p}\big\|H_{\Phi,\vec A}^{\vec b}(\vec f)\chi_k\big\|^p_{L^{q(\cdot)}_{\omega}}\Big)^{\frac{1}{p}},\nonumber
\\
&&E_2 = 2^{-k_0\lambda}\Big(\sum\limits_{k=-\infty}^{-1}2^{k\alpha^*(0)p}\big\|H_{\Phi,\vec A}^{\vec b}(\vec f)\chi_k\big\|^p_{L^{q(\cdot)}_{\omega}}\Big)^{\frac{1}{p}},\nonumber
\\
&&E_3=2^{-k_0\lambda}\Big(\sum\limits_{k=0}^{k_0}2^{k\alpha^*_\infty p}\big\|H_{\Phi,\vec A}^{\vec b}(\vec f)\chi_k\big\|^p_{L^{q(\cdot)}_{\omega}}\Big)^{\frac{1}{p}}\nonumber.
\end{eqnarray}
Now, we need to estimate the upper bounds for $E_1, E_2$ and $E_3$. Note that, 
using (\ref{Hvec4}), $E_1$ is dominated by
\begin{eqnarray}\label{estimateE1}
E_1&\lesssim& \mathcal C_1.\mathcal B_{\rm Lip}.\mathcal F.2^{-k_0\lambda}\Big(\sum\limits_{k=-\infty}^{k_0}2^{k(\alpha^*(0)+\beta+\sum\limits_{i=1}^m (\gamma_i+n)\sigma_i)p}\prod_{i=1}^{m}\big(2^{k(\lambda_i -\alpha_i(0))p}+2^{k(\lambda_i -\alpha_{i\infty})p}\big)\Big)^{\frac{1}{p}}\nonumber
\\
&:=&  \mathcal C_1.\mathcal B_{\rm Lip}.\mathcal F.\mathcal T_0.
\end{eqnarray}
In view of $\alpha_*$, we have
\[
\mathcal T_0 = 2^{-k_0\lambda}.\Big(\sum\limits_{k=-\infty}^{k_0}2^{k(\sum\limits_{i=1}^m\alpha_i(0)+\sum\limits_{i=1}^m (\gamma_i+n)(\sigma_i-\frac{1}{r_i(0)})p}\prod\limits_{i=1}^{m}\big(2^{k(\lambda_i -\alpha_i(0))p}+2^{k(\lambda_i -\alpha_{i\infty})p}\big)\Big)^{\frac{1}{p}}.
\]
Note that, by defining $\sigma_i$ and (\ref{DKgammari}), it is clear to see that 
\begin{equation}\label{sigma}
(\gamma_i+n)(\sigma_i-\frac{1}{r_i(0)})=0,\,\textit{\rm for all}\, i=1,...,m.
\end{equation}
So, we get
\begin{eqnarray}
\mathcal T_0 &=&2^{-k_0\lambda}.\Big(\sum\limits_{k=-\infty}^{k_0}2^{k\big(\sum\limits_{i=1}^m\alpha_i(0)\big)p}\prod\limits_{i=1}^{m}\big(2^{k(\lambda_i -\alpha_i(0))p}+2^{k(\lambda_i -\alpha_{i\infty})p}\big)\Big)^{\frac{1}{p}}\nonumber
\\
&=& 2^{-k_0\lambda}\Big(\sum\limits_{k=-\infty}^{k_0}\prod_{i=1}^{m}\big(2^{k\lambda_i p}+2^{k(\lambda_i -\alpha_{i\infty}+\alpha_i(0))p}\big)\Big)^{\frac{1}{p}}\nonumber
\\
&\lesssim& \Big(\prod_{i=1}^{m} 2^{-k_0\lambda_i p}\Big\{  \sum\limits_{k=-\infty}^{k_0} 2^{k\lambda_i p}+ \sum\limits_{k=-\infty}^{k_0}2^{k(\lambda_i -\alpha_{i\infty}+\alpha_i(0))p}\Big\}\Big)^{\frac{1}{p}}.\nonumber
\end{eqnarray}
Because of assuming that $\lambda_i >0$, for all $i=1,...,m$ and (\ref{DKalpha}), we obtain
\begin{eqnarray}
\mathcal T_0 &\lesssim& \Big(\prod_{i=1}^{m}2^{-k_0\lambda_i p}\Big\{\dfrac{2^{k_0\lambda_i p}}{1-2^{-\lambda_ip}}+ \dfrac{2^{k_0(\lambda_i -\alpha_{i\infty}+\alpha_i(0))p}}{1- 2^{-(\lambda_i -\alpha_{i\infty}+\alpha_i(0))p}}\Big\}\Big)^{\frac{1}{p}}\nonumber
\\
&\lesssim &\prod_{i=1}^{m}\Big\{\dfrac{1}{1-2^{-\lambda_ip}}+ \dfrac{2^{k_0(-\alpha_{i\infty}+\alpha_i(0))}}{1- 2^{-(\lambda_i -\alpha_{i\infty}+\alpha_i(0))p}}\Big\}\lesssim \prod\limits_{i=1}^{m}\Big(1+ 2^{{k_0}\big(\alpha_i(0)-\alpha_{i\infty}\big)}\Big).\nonumber
\end{eqnarray}
Consequently, from (\ref{estimateE1}), we conclude 
\begin{equation}\label{estimateE1*}
E_1\lesssim \mathcal C_1.\mathcal B_{\rm Lip}.\mathcal F.\prod\limits_{i=1}^{m}\Big(1+ 2^{{k_0}\big(\alpha_i(0)-\alpha_{i\infty}\big)}\Big).
\end{equation}
A similar agrument as $E_1$, we also get
\begin{equation}\label{estimateE2*}
E_2\lesssim \mathcal C_1.\mathcal B_{\rm Lip}.\mathcal F.2^{-k_0\lambda}.
\end{equation}
For $i=1,...,m$, we define
\[
{L_i} = \left\{ \begin{gathered}
{2^{{k_0}({\alpha _{i\infty} } - \alpha_i(0))}} + {\left| {{2^{\lambda_i p}} - 1} \right|^{ - \frac{1}{p}}} + {2^{ - {k_0}\lambda_i }},\,\textit{\rm if}\,\lambda_i  + {\alpha_{i\infty}} - \alpha_i(0) \ne 0, \hfill \\
{2^{ - {k_0}\lambda_i }}{({k_0} + 1)^{\frac{1}{p}}} + {\left| {{2^{\lambda_i p}} - 1} \right|^{ - \frac{1}{p}}},\, \rm otherwise .\hfill \\ 
\end{gathered}  \right.
\]
Thus, we see that
\begin{equation}\label{estimateE3}
E_3 \lesssim \mathcal C_1.\mathcal B_{\rm Lip}.\mathcal F.\mathcal T_\infty,
\end{equation}
where 
\[\mathcal T_\infty= 2^{-k_0\lambda}\Big(\sum\limits_{k=0}^{k_0}2^{k(\sum\limits_{i=1}^m\alpha_{i\infty}+\sum\limits_{i=1}^m (\gamma_i+n)(\sigma_i-\frac{1}{r_{i\infty}})p}\prod\limits_{i=1}^{m}\big(2^{k(\lambda_i -\alpha_i(0))p}+2^{k(\lambda_i -\alpha_{i\infty})p}\big)\Big)^{\frac{1}{p}}.
\]
Remark that, by defining $\sigma_i$ and (\ref{DKgammari}), we deduce 
\begin{equation}\label{sigma1}
(\gamma_i+n)(\sigma_i-\frac{1}{r_{i\infty}})=0,\,\textit{\rm for all}\, i=1,...,m.
\end{equation}
Thus, by estimating in the same way as  $\mathcal T_0$, we also have
\begin{eqnarray}
\mathcal T_\infty &\lesssim& \prod_{i=1}^{m} 2^{-k_0\lambda_i}\Big( \sum\limits_{k=0}^{k_0} 2^{k\lambda_i p}+ \sum\limits_{k=0}^{k_0}2^{k(\lambda_i +\alpha_{i\infty}-\alpha_i(0))p}\Big)^{\frac{1}{p}}:= \prod\limits_{i=1}^{m} \mathcal T_{i,\infty}.
\nonumber
\end{eqnarray}
In the case $\lambda_i+\alpha_{i\infty}-\alpha_i(0)= 0$, $\mathcal T_{i,\infty}$ is dominated by
\[
\mathcal T_{i,\infty} \leq 2^{-k_0\lambda_i} \Big( \dfrac{2^{k_0\lambda_i p}-1}{2^{\lambda_i.p}-1}+(k_0+1)\Big)^{\frac{1}{p}}\lesssim  {2^{ - {k_0}\lambda_i }}{({k_0} + 1)^{\frac{1}{p}}} + {\left| {{2^{\lambda_i p}} - 1} \right|^{ - \frac{1}{p}}}.
\]
Otherwise, we get 
\begin{eqnarray}
\mathcal T_{i,\infty}&\leq& 2^{-k_0\lambda_i}\Big( \dfrac{2^{k_0\lambda_i p}-1}{2^{\lambda_i.p}-1}+\dfrac{2^{k_0(\lambda_i +\alpha_{i\infty}-\alpha_i(0))p}-1}{2^{(\lambda_i +\alpha_{i\infty}-\alpha_i(0))p}-1}\Big)^{\frac{1}{p}} \nonumber
\\
&\lesssim & {2^{{k_0}({\alpha _{i\infty} } - \alpha_i(0))}} + {\left| {{2^{\lambda_i p}} - 1} \right|^{ - 1/p}} + {2^{ - {k_0}\lambda_i }}\nonumber.
\end{eqnarray}
which implies $\mathcal T_\infty\lesssim \prod\limits_{i=1}^m L_i.$ From this, by (\ref{estimateE3}), we obtain
\begin{equation}\label{estimateE3*}
E_3\lesssim \mathcal C_1.\mathcal B_{\rm Lip}.\mathcal F.\prod\limits_{i=1}^m L_i.
\end{equation}
By (\ref{Hvec5}), (\ref{estimateE1*}), (\ref{estimateE2*}) and (\ref{estimateE3*}), the proof of Theorem \ref{TheoremMorreyHerzLip} is finished.
\vskip 5pt
Next, let us give the proof for Theorem \ref{TheoremHerzLip}.
 From Proposition 3.8 in \cite{Almeida2012}, it is easy to see that
\begin{eqnarray}\label{H0H1}
{\big\| {{H_{\Phi ,\vec A }^{\vec b}}\big( {\vec f } \big)} \big\|_{{\mathop K\limits^.}_{q(\cdot),\omega}^{{\alpha^*(\cdot)}, p}}}&\lesssim & \Big(\sum\limits_{k=-\infty}^{-1}{2^{k\alpha^*(0)p}\big\|{{H_{\Phi ,\vec A }}\big( {\vec f } \big)}\chi_k\big\|_{L^{q(\cdot)}_{\omega}}^p}\Big)^{\frac{1}{p}}\nonumber
\\
&&\,\,\,\,\,\,\,\,\,\,\,\,\,\,+ \Big(\sum\limits_{k=0}^{\infty}{2^{k\alpha^*_\infty p}\big\|{{H_{\Phi ,\vec A }}\big( {\vec f } \big)}\chi_k\big\|_{L^{q(\cdot)}_{\omega}}^p}\Big)^{\frac{1}{p}}\nonumber
\\\
&:=&\mathcal H_0+\mathcal H_1.
\end{eqnarray} 
Next, we need to estimate  the upper bound of $\mathcal H_0$ and $\mathcal H_1$. In view of (\ref{Hvec2}) and (\ref{sigma}), by using the Minkowski inequality,  we find
\begin{eqnarray}\label{H0L1}
\mathcal H_0 &\lesssim& \mathcal B_{\rm Lip}{\int\limits_{{\mathbb R^n}} {\frac{\Phi(t)}{|t|^n}\prod\limits_{i=1}^{m} c_{A_i,q_i,\gamma_i}(t)\big\|1\big\|_{L^{\vartheta_i(t,\cdot)}}}}\big\|I_n-A_i(t)\big\|^{\beta_i}\times
\\
&&\,\,\,\,\,\,\,\,\,\,\,\,\,\,\,\,\times\Big\{ {{{\sum\limits_{k =  - \infty }^{{-1}} {{2^{k(\sum\limits_{i=1}^m\alpha_i(0))p}}\prod\limits_{i = 1}^m {\Big( {\sum\limits_{r = {\Theta _n^*} - 1}^0 {\left\| {{f_i}{\chi _{k + {\ell_i} + r}}} \right\|_{L^{\zeta q_i(\cdot)}_{\omega_i}}} } \Big)^p} } }}} \Big\} ^{\frac{1}{p}}dt.
\nonumber
\end{eqnarray}
Using (\ref{DKpip}) and the H\"{o}lder inequality, it follows that
\begin{eqnarray}\label{Holderlr}
&&{\Big\{ {{{\sum\limits_{k = - \infty }^{{-1}} {{2^{k(\sum\limits_{i=1}^m\alpha_i(0))p}}\prod\limits_{i = 1}^m {\Big( {\sum\limits_{r = {\Theta _n^*} - 1}^0 {\left\| {{f_i}{\chi _{k + {\ell_i} + r}}} \right\|_{L^{\zeta q_i(\cdot)}_{\omega_i}}} } \Big)^p} } }}} \Big\}^{\frac{1}{p}}}
\\
&&\,\,\,\,\,\,\,\,\leq \prod\limits_{i = 1}^m {{{\Big\{ {{{\sum\limits_{k =  - \infty }^{{-1}} {{2^{k{\alpha_i(0)}{p_i}}}\Big( {\sum\limits_{r = {\Theta_n^*} - 1}^0 {\left\| {{f_i}{\chi _{k + {\ell_i} + r}}} \right\|_{L^{\zeta q_i(\cdot)}_{\omega_i}}} } \Big)^{p_i}} }}} \Big\}}^{\frac{1}{{p_i}}}}}.\nonumber
\end{eqnarray}
By $p_i\geq 1$, for all $i=1,...,m$, we have
\[
{\Big( {\sum\limits_{r = {\Theta _n^*} - 1}^0 {\left\| {{f_i}{\chi _{k + {\ell_i} + r}}} \right\|_{L^{\zeta q_i(\cdot)}_{\omega_i}}^{}} } \Big)^{{p_i}}} \leq {\left( {2 - {\Theta _n^*}} \right)^{{p_i} - 1}}\sum\limits_{r = {\Theta _n^*} - 1}^0 {\left\| {{f_i}{\chi _{k + {\ell_i} + r}}} \right\|_{L^{\zeta q_i(\cdot)}_{\omega_i}}^{{p_i}}}.
\]
Thus, combining (\ref{H0L1}) and (\ref{Holderlr}), we deduce
\begin{equation}\label{H0L2}
\mathcal H_0 \lesssim \mathcal B_{\rm Lip}.\int\limits_{{\mathbb R^n}} {{(2 - {\Theta_n^*})^{m - \frac{1}{p}}}\frac{\Phi(t)}{|t|^n}\prod\limits_{i=1}^{m}c_{A_i,q_i,\gamma_i}(t)\big\|1\big\|_{L^{\vartheta_i(t,\cdot)}}\big\|I_n-A_i(t)\big\|^{\beta_i}{\mathcal{H}_{0,i}}}dt.
\end{equation}
Here $\mathcal{H}_{0,i} = {\sum\limits_{r = {\Theta _n^*} - 1}^0 {{{\Big( {\sum\limits_{k =  - \infty }^{{-1}} {{2^{k{\alpha _i(0)}{p_i}}}\left\| {{f_i}{\chi _{k + {\ell_i} + r}}} \right\|_{L^{\zeta q_i(\cdot)}_{\omega_i}}^{{p_i}}} } \Big)}^{^{\frac{1}{{p_i}}}}}} }$ for all $i=1,2,...,m$. 
\\
Hence, we estimate
\begin{eqnarray}\label{Hi}
{\mathcal{H}_{0,i}} &=& {\sum\limits_{r = {\Theta _n^*} - 1}^0 {{{\Big( {\sum\limits_{t = - \infty }^{-1+\ell+r} {{2^{(t-\ell_i-r){\alpha _i(0)}{p_i}}}\left\| {{f_i}{\chi _{t}}} \right\|_{L^{\zeta q_i(\cdot)}_{\omega_i}}^{{p_i}}} } \Big)}^{^{\frac{1}{{p_i}}}}}} }
\nonumber
\\
 &\lesssim&{\sum\limits_{r = {\Theta _n^*} - 1}^0 { 2^{-(\ell_i+r)\alpha_i(0)}{{\Big( {\sum\limits_{t =  - \infty }^{\infty} {{2^{t{\alpha _i(0)}{p_i}}}\left\| {{f_i}{\chi _{t}}} \right\|_{L^{\zeta q_i(\cdot)}_{\omega_i}}^{{p_i}}} } \Big)}^{^{\frac{1}{{p_i}}}}}}}.
\nonumber
\end{eqnarray}
By $\alpha_i(0)=\alpha_{i\infty}$ and Proposition 3.8 in \cite{Almeida2012}, we get
\begin{equation}\label{H0i}
{\mathcal{H}_{0,i}}\lesssim {\sum\limits_{r = {\Theta _n^*} - 1}^0 { 2^{-(\ell_i+r)\alpha_i(0)}}}{{{\left\| {{f_i}} \right\|}_{{\mathop K\limits^.}_{\zeta q_i(\cdot),\omega_i}^{{\alpha_i(\cdot)}, p_i}}}}=2^{-\ell_i\alpha_i(0)}.{\sum\limits_{r = {\Theta _n^*} - 1}^0 { 2^{-r\alpha_i(0)}}}{{{\left\| {{f_i}} \right\|}_{{\mathop K\limits^.}_{\zeta q_i(\cdot),\omega_i}^{{\alpha_i(\cdot)}, p_i}}}}.
\end{equation}
Since $2^{\ell_i-1}<\left\|A_i(t)\right\|\leq 2^{\ell_i}$, we deduce that ${2^{^{{-\ell_i}{\alpha _i(0)}}}} \lesssim {\left\| {{A_i}(t)} \right\|^{- {\alpha_i(0)}}}$. Hence, by (\ref{H0i}), we have
\[
\mathcal H_{0,i}\lesssim \phi_{A_i,0}(t).{\left\| {{f_i}} \right\|}_{{\mathop K\limits^.}_{\zeta q_i(\cdot),\omega_i}^{{\alpha_i(\cdot)}, p_i}}.
\]
As above, by (\ref{H0L2}), we make
\[
\mathcal H_0\lesssim \mathcal C_2.\mathcal B_{\rm Lip}.\prod\limits_{i = 1}^m {{{\left\| {{f_i}} \right\|}_{{\mathop K\limits^.}_{\zeta q_i(\cdot),\omega_i}^{{\alpha_i(\cdot)}, p_i}}}}.
\]
By estimating as $\mathcal H_0$, we also make 
\[
\mathcal H_1\lesssim \mathcal C_2.\mathcal B_{\rm Lip}.\prod\limits_{i = 1}^m {{{\left\| {{f_i}} \right\|}_{{\mathop K\limits^.}_{\zeta q_i(\cdot),\omega_i}^{{\alpha_i(\cdot)}, p_i}}}}.
\]
There, by (\ref{H0H1}), we finishes desired conclusion.
\\
\subsection{Proofs of Theorem \ref{TheoremMorreyHerzCMO} and Theorem \ref{TheoremHerCMO}}
Applying the Minkowski inequality and the H\"{o}lder inequality  for variable Lebesgue spaces, we get
\[
\big\|H_{\Phi,\vec{A}}^{\vec b}(\vec{f})\chi_k\big\|_{L^{q(\cdot)}_{\omega}}\lesssim \int\limits_{\mathbb R^n}{\frac{\Phi(t)}{|t|^n}}\prod_{i=1}^{m}\big\|\big(b_i(\cdot)-b_i(A_i(t)\cdot)\big)\big\|_{L^{r_i}(\omega_i,B_k)}\big\|f_i(A_i(t).)\chi_k\big\|_{L^{q_i(\cdot)}_{\omega_i}}dt.
\]
By(\ref{fAchikLq}), we deduce 
\begin{eqnarray}\label{Hvec2CMO}
\big\|H_{\Phi,\vec A}^{\vec b}(\vec f)\chi_k\big\|_{L^{q(\cdot)}_{\omega}}&\lesssim &\int\limits_{\mathbb R^n}{\frac{\Phi(t)}{|t|^n}\prod\limits_{i=1}^{m}}c_{A_i,q_i,\gamma_i}(t)\big\|b_i(\cdot)-b_i(A_i(t)\cdot)\big\|_{L^{r_i}(\omega_i,B_k)}\times
\nonumber
\\
&&\,\,\times\big\|1\big\|_{L^{\vartheta_i(t,\cdot)}}\prod\limits_{i=1}^{m}\sum\limits_{r=\Theta_n^*-1}^{0} \big\|f_i\chi_{k+\ell_i+r}\big\|_{L^{\zeta q_i(\cdot)}_{\omega_i}}dt.
\end{eqnarray}
On the other hand, we need to prove that
\begin{equation}\label{BDTCMO}
\big\|b_i(\cdot)-b_i(A_i(t)\cdot)\big\|_{L^{r_i}(\omega_i, B_k)}\lesssim 2^{\frac{k(\gamma_i+n)}{r_i}}\Big(1+\psi_{A_i,\gamma_i}^{\frac{1}{r_i}}|s_i(t)|^{\frac{\gamma_i+n}{r_i}}+\varphi_{A_i}(t)\Big)\big\|b_i\big\|_{{\mathop {CMO}\limits^{.}}^{r_i}(\omega_i)}.
\end{equation}
In fact, we put $ a_{1,i}(\cdot)= b_i(\cdot)-b_{i,\omega_i,B_k},$ $ a_{2,i}(\cdot)= b_i(A_i(t).)-b_{i,\omega_i,A_i(t)B_k}$ and $ a_{3,i}(\cdot)= b_{i,\omega_i,B_k}-b_{i,\omega_i,A_i(t)B_k}.$ Here
\[
 b_{i,\omega_i,U}=\frac{1}{\omega_i(U)}\int\limits_{U}{b_i(x)\omega_i(x)dx}.
 \]
Then, we have
\begin{equation}\label{BDTphantichCMO}
\big\|b_i(\cdot)-b_i(A_i(t)\cdot)\big\|_{L^{r_i}{(\omega_i,B_k)}}\leq \big\|a_{1,i}\big\|_{L^{r_i}(\omega_i,B_k)}+\big\|a_{2,i}\big\|_{L^{r_i}(\omega_i,B_k)}+\big\|a_{3,i}\big\|_{L^{r_i}(\omega_i,B_k)}.
\end{equation}
From defining the space ${\mathop{CMO}\limits^{.}}^{r_i}(\omega_i)$, we immediately have
\begin{equation}\label{a1i}
\big\|a_{1,i}\big\|_{L^{r_i}(\omega_i,B_k)}\leq \big(\omega_i(B_k)\big)^{\frac{1}{r_i}}.\big\|b_i\big\|_{{\mathop{ CMO}\limits^{.}}^{r_i}(\omega_i)}\lesssim 2^{\frac{k(\gamma_i+n)}{r_i}}.\big\|b_i\big\|_{{\mathop{CMO}\limits^{.}}^{r_i}(\omega_i)}.
\end{equation}
By making the formula for change of variables , we obtain
\begin{eqnarray}
\big\|a_{2,i}\big\|^{r_i}_{L^{r_i}(\omega_i,B_k)}&=&\int\limits_{B_k}{|b_i(A_i(t)x)-b_{i,\omega_i,A_i(t)B_k}|^{r_i}\omega_i(x)dx}
\nonumber
\\
&\leq & \psi_{A_i,\gamma_i}(t).\int\limits_{A_i(t)B_k}{|b_i(z)-b_{i,\omega_i,A_i(t)B_k}|^{r_i}\omega_i(z)dx}.
\nonumber
\end{eqnarray}
Because of assuming $A_i(t)=s_i(t).a_i(t)$, we deduce
\begin{equation}\label{a2i}
\big\|a_{2,i}\big\|_{L^{r_i}(\omega_i,B_k)}\lesssim \psi_{A_i,\gamma_i}^{\frac{1}{r_i}}|s_i(t)|^{\frac{\gamma_i+n}{r_i}}.2^{\frac{k(\gamma_i+n)}{r_i}}\big\|b_i\big\|_{{\mathop{CMO}\limits^{.}}^{r_i}(\omega_i)}.
\end{equation}
Next, we observe that
\begin{equation}\label{a3i}
\big\|a_{3,i}\big\|_{L^{r_i}(\omega_i, B_k)}\leq(\omega_i(B_k))^{\frac{1}{r_i}}\big|b_{i,\omega_i,B_k}-b_{i,\omega_i,A_i(t)B_k}\big|.
\end{equation}
By having $s_i(t)\neq 0$, there exists an integer number $\theta_i = \theta_i(t)$ satisfying $2^{\theta_i-1}<|s_i(t)|\leq 2^{\theta_i}$. Thus, we define 
\[
\sigma(\theta_i) = \left\{ \begin{array}{l}
\theta_i-1,\,\textit{\rm if}\,\theta_i\geq 1,
\\
\theta_i,\,\,\,\,\,\,\,\,\,\,\,\textit{\rm otherwise},
\end{array} \right.
\]
and 
\[
S(\theta_i) = \left\{ \begin{array}{l}
\big\{j\in\mathbb Z: 1\leq j\leq \theta_i-1\big\},\,\textit{\rm if}\,\theta_i\geq 1,
\\
\\
\big\{j\in\mathbb Z: \theta_i+1\leq j\leq 0\big\},\textit{\rm otherwise}.
\end{array} \right.
\]
At this point, we give the estimation as below
\begin{equation}\label{bdtBMOold}
\big|b_{i,\omega_i,B_k}-b_{i,\omega_i,A_i(t)B_k}\big|\leq \sum\limits_{j\in S(\theta_i)}\big|b_{i,\omega_i,2^{j-1}B_k}-b_{i,\omega_i,2^jB_k}\big|+\big|b_{i,\omega_i,2^{\sigma(\theta_i)}B_k}-b_{i,\omega_i,A_i(t)B_k}\big|.
\end{equation}
When $S(\theta_i)$ is empty set, we should understand that
\[
 \sum\limits_{j\in S(\theta_i)}\big|b_{i,\omega_i,2^{j-1}B_k}-b_{i,\omega_i,2^jB_k}\big|:=0.
\]
It is not difficult to show that
\[
\big|b_{i,\omega_i,2^{j-1}B_k}-b_{i,\omega_i,2^jB_k}\big|\lesssim \big\|b_i\big\|_{{\mathop{CMO}\limits^{.}}^{r_i}(\omega_i)}.
\]
In the case $\theta_i\geq 1$, by defining $\sigma$, it follows that
\begin{eqnarray}
\big|b_{i,\omega_i,2^{\sigma(\theta_i)}B_k}-b_{i,\omega_i,A_i(t)B_k}\big| &\leq & \frac{1}{\omega_i(2^{\theta_i-1}B_k)}\int\limits_{2^{\theta_i-1}B_k}{\big|b_i(x)-b_{i,\omega_i,A_i(t)B_k}\big|\omega_i(x)dt}\nonumber
\\
&\leq&\frac{\big(\omega_i(A_i(t)B_k)\big)^{\frac{1}{r_i'}}}{\omega_i(2^{\theta_i-1}B_k)}\Big(\int\limits_{A_i(t)B_k}{\big|b_i(x)-b_{i,\omega_i,A_i(t)B_k}\big|^{r_i}\omega_i(x)dt}\Big)^{\frac{1}{r_i}}
\nonumber
\\
&\leq&\frac{\big(\omega_i(A_i(t)B_k)\big)}{\omega_i(2^{\theta_i-1}B_k)}\big\|b_i\big\|_{{\mathop{CMO}\limits^{.}}^{r_i}(\omega_i)}\nonumber.
\end{eqnarray}
Note that $A_i(t)=s_i(t)a_i(t)$. So, we compute
\begin{equation}\label{tile}
\frac{\big(\omega_i(A_i(t)B_k)\big)}{\omega_i(2^{\theta_i-1}B_k)}\lesssim \frac{\big(|s_i(t)|2^k\big)^{n+\gamma_i}}{(2^{\theta_i-1}2^k)^{n+\gamma_i}}\leq \frac{\big(2^{\theta_i}2^k\big)^{n+\gamma_i}}{(2^{\theta_i-1}2^k)^{n+\gamma_i}}\lesssim 1.
\end{equation}
Consequently, we have
\[
\big|b_{i,\omega_i,2^{\sigma(\theta_i)}B_k}-b_{i,\omega_i,A_i(t)B_k}\big|\lesssim \big\|b_i\big\|_{{\mathop{CMO}\limits^{.}}^{r_i}(\omega_i)}. 
\]
Otherwise, for $\theta_i\leq 0$, by estimating as (\ref{tile}), we deduce
\begin{eqnarray}
\big|b_{i,\omega_i,2^{\sigma(\theta_i)}B_k}-b_{i,\omega_i,A_i(t)B_k}\big| &\leq & \frac{1}{\omega_i(A_i(t)B_k)}\int\limits_{A_i(t)B_k}{\big|b_i(x)-b_{i,\omega_i,2^{\theta_i}a_i(t)B_k}\big|\omega_i(x)dt}\nonumber
\\
&\leq&\frac{\big(\omega_i(2^{\theta_i}a_i(t)B_k)\big)^{\frac{1}{r_i'}}}{\omega_i(A_i(t)B_k)}\Big(\int\limits_{2^{\theta_i}a_i(t)B_k}{\big|b_i(x)-b_{i,\omega_i,2^{\theta_i}a_i(t)B_k}\big|^{r_i}\omega_i(x)dt}\Big)^{\frac{1}{r_i}}
\nonumber
\\
&\leq&\frac{\big(\omega_i(2^{\theta_i}a_i(t)B_k)\big)}{\omega_i(A_i(t)B_k)}\big\|b_i\big\|_{{\mathop{CMO}\limits^{.}}^{r_i}(\omega_i)}\lesssim \big\|b_i\big\|_{{\mathop{CMO}\limits^{.}}^{r_i}(\omega_i)}.
\nonumber
\end{eqnarray}
Because of  $2^{\theta_i-1}<|s_i(t)|\leq 2^{\theta_i}$, we have
\[
|\theta_i|+1 \lesssim \left\{ \begin{array}{l}
{\rm log}(4|s_i(t)|),\,\,\textit{\rm if}\,\,\theta_i\geq 0,
\\
\\
{\rm log}\frac{2}{|s_i(t)|},\,\,\textit{\rm otherwise}.
\end{array} \right.
\lesssim \varphi_{A_i}(t).
\]
Therefore, by having (\ref{bdtBMOold}), it follows that
\[
\big|b_{i,\omega_i,B_k}-b_{i,\omega_i,A_i(t)B_k}\big|\lesssim (|\theta_i|+1)\big\|b_i\big\|_{{\mathop{CMO}\limits^{.}}^{r_i}(\omega_i)}\lesssim \varphi_{A_i}(t).\big\|b_i\big\|_{{\mathop{CMO}\limits^{.}}^{r_i}(\omega_i)}.
\]
As above, by (\ref{a3i}), we get
\[
\big\|a_{3,i}\big\|_{L^{r_i}(\omega_i,B_k)}\lesssim 2^{\frac{k(n+\gamma_i)}{r_i}}\varphi_{A_i}(t).\big\|b_i\big\|_{{\mathop{CMO}\limits^{.}}^{r_i}(\omega_i)}.
\]
From this, by (\ref{a1i}), (\ref{a2i}), we finish the proof of the inequality (\ref{BDTCMO}).
\\
Using (\ref{Hvec2CMO}) and (\ref{BDTCMO}), we have
\begin{eqnarray}\label{HMHCMO}
\big\|H_{\Phi,\vec A}^{\vec b}(\vec f)\chi_k\big\|_{L^{q(\cdot)}_{\omega}}&\lesssim &\mathcal B_{CMO, \vec{\omega}}. 2^{k(\sum\limits_{i=1}^m \frac{\gamma_i+n}{r_i})}\Big(\int\limits_{\mathbb R^n}{\frac{\Phi(t)}{|t|^n}\prod\limits_{i=1}^{m}}c_{A_i,q_i,\gamma_i}(t)\big\|1\big\|_{L^{\vartheta_i(t,\cdot)}}\times
\nonumber
\\
&\times&\Big(1+\psi_{A_i,\gamma_i}^{\frac{1}{r_i}}|s_i(t)|^{\frac{\gamma_i+n}{r_i}}+\varphi_{A_i}(t)\Big)\prod\limits_{i=1}^{m}\sum\limits_{r=\Theta_n^*-1}^{0} \big\|f_i\chi_{k+\ell_i+r}\big\|_{L^{\zeta q_i(\cdot)}_{\omega_i}}dt\Big).
\nonumber
\\
\end{eqnarray}
At this point, by making Lemma \ref{lemmaVE} in Section 2 again, we have a similar results to (\ref{Hvec4}) as follow
\[
\big\|H_{\Phi,\vec A}^{\vec b}(\vec f)\chi_k\big\|_{L^{q(\cdot)}_{\omega}}\lesssim\mathcal C_4.\mathcal B_{CMO,\vec{\omega}}.\mathcal F. 2^{k(\sum\limits_{i=1}^m \frac{\gamma_i+n}{r_i})}\prod_{i=1}^{m}\big(2^{k(\lambda_i -\alpha_i(0))}+2^{k(\lambda_i -\alpha_{i\infty})}\big).
\]
By using Proposition 2.5 in \cite{LZ2014} again, we get
\begin{equation}\label{HMHCM01}
\big\| H_{\Phi,\vec A}^{\vec b}(\vec f)\big\|_{M{\mathop K\limits^.}^{\alpha^{**}(\cdot), \lambda}_{p,q(\cdot),\omega}} \lesssim \max\big\{
\sup\limits_{k_0<0,k_0\in\mathbb Z} {\widetilde E}_1, \sup\limits_{k_0\geq 0,k_0\in\mathbb Z}({\widetilde E}_2+ {\widetilde E}_3)
\big\},
\end{equation}
where 
\begin{eqnarray}
&&\widetilde E_1=2^{-k_0\lambda}\Big(\sum\limits_{k=-\infty}^{k_0}2^{k\alpha^{**}(0)p}\big\|H_{\Phi,\vec A}(\vec f)\chi_k\big\|^p_{L^{q(\cdot)}_{\omega}}\Big)^{\frac{1}{p}},\nonumber
\\
&&\widetilde E_2 = 2^{-k_0\lambda}\Big(\sum\limits_{k=-\infty}^{-1}2^{k\alpha^{**}(0)p}\big\|H_{\Phi,\vec A}(\vec f)\chi_k\big\|^p_{L^{q(\cdot)}_{\omega}}\Big)^{\frac{1}{p}},\nonumber
\\
&&\widetilde E_3=2^{-k_0\lambda}\Big(\sum\limits_{k=0}^{k_0}2^{k\alpha_\infty^{**} p}\big\|H_{\Phi,\vec A}(\vec f)\chi_k\big\|^p_{L^{q(\cdot)}_{\omega}}\Big)^{\frac{1}{p}}\nonumber.
\end{eqnarray}
In view of (\ref{HMHCM01}), by defining $\alpha^{**}$ and estimating as (\ref{estimateE1*}), (\ref{estimateE2*}), (\ref{estimateE3*}), we also have
\begin{eqnarray}\label{estimateE123CMO}
\widetilde E_1 &\lesssim& \mathcal C_4.\mathcal B_{CMO,\vec{\omega}}.\mathcal F.\prod\limits_{i=1}^{m}\Big(1+ 2^{{k_0}\big(\alpha_i(0)-\alpha_{i\infty}\big)}\Big),
\nonumber
\\
\widetilde E_2 &\lesssim& \mathcal C_4.\mathcal B_{CMO,\vec{\omega}}.\mathcal F.2^{-k_0\lambda},
\nonumber
\\
\widetilde E_3 &\lesssim& \mathcal C_4.\mathcal B_{CMO,\vec{\omega}}.\mathcal F.\prod\limits_{i=1}^m L_i.
\nonumber
\end{eqnarray}
Therefore, the proof of Theorem \ref{TheoremMorreyHerzCMO} is completed.
\vskip 5pt
Now, let us give the proof for Theorem \ref{TheoremHerCMO}. By making Proposition 3.8 in \cite{Almeida2012} again, we obtain
\begin{eqnarray}\label{H0H1CMO}
{\big\| {{H_{\Phi ,\vec A }^{\vec b}}\big( {\vec f } \big)} \big\|_{{\mathop K\limits^.}_{q(\cdot),\omega}^{{\alpha^{**}(\cdot)}, p}}}&\lesssim & \Big(\sum\limits_{k=-\infty}^{-1}{2^{k\alpha^{**}(0)p}\big\|{{H_{\Phi ,\vec A }^{\vec b}}\big( {\vec f } \big)}\chi_k\big\|_{L^{q(\cdot)}_{\omega}}^p}\Big)^{\frac{1}{p}}\nonumber
\\
&&\,\,\,\,\,\,\,\,\,\,\,\,\,\,+ \Big(\sum\limits_{k=0}^{\infty}{2^{k\alpha^{**}_\infty p}\big\|{{H_{\Phi ,\vec A }^{\vec b}}\big( {\vec f } \big)}\chi_k\big\|_{L^{q(\cdot)}_{\omega}}^p}\Big)^{\frac{1}{p}}\nonumber
\\\
&:=&\mathcal G_0+\mathcal G_1.
\end{eqnarray} 
Using the Minkowski inequality, by employing (\ref{HMHCMO}), we find
\begin{eqnarray}\label{H0L1CMO}
\mathcal G_0 &\lesssim& \mathcal B_{CMO,\vec{\omega}}{\int\limits_{{\mathbb R^n}} {\frac{\Phi(t)}{|t|^n}\prod\limits_{i=1}^{m}c_{A_i,q_i,\gamma_i}(t)\big\|1\big\|_{L^{\vartheta_i(t,\cdot)}}}}\phi_{A_i,0}(t)\big(1+\psi_{A_i,\gamma_i}^{\frac{1}{r_i}}|s_i(t)|^{\frac{\gamma_i+n}{r_i}}+\varphi_{A_i}(t)\big)\times
\nonumber
\\
&&\,\,\,\,\,\,\,\,\,\,\,\,\,\,\,\,\times\Big\{ {{{\sum\limits_{k =  - \infty }^{{-1}} {{2^{k(\sum\limits_{i=1}^m\alpha_i(0))p}}\prod\limits_{i = 1}^m {\Big( {\sum\limits_{r = {\Theta _n^*} - 1}^0 {\left\| {{f_i}{\chi _{k + {\ell_i} + r}}} \right\|_{L^{\zeta q_i(\cdot)}_{\omega_i}}} } \Big)} } }^p}} \Big\} ^{\frac{1}{p}}dt,
\nonumber
\end{eqnarray}
and
\begin{eqnarray}\label{H0L1CMO}
\mathcal G_1 &\lesssim& \mathcal B_{CMO,\vec{\omega}}{\int\limits_{{\mathbb R^n}} {\frac{\Phi(t)}{|t|^n}\prod\limits_{i=1}^{m}c_{A_i,q_i,\gamma_i}(t).\big\|1\big\|_{L^{\vartheta_i(t,\cdot)}}}}\phi_{A_i,0}(t)\big(1+\psi_{A_i,\gamma_i}^{\frac{1}{r_i}}|s_i(t)|^{\frac{\gamma_i+n}{r_i}}+\varphi_{A_i}(t)\big)\times
\nonumber
\\
&&\,\,\,\,\,\,\,\,\,\,\,\,\,\,\,\,\times\Big\{ {{{\sum\limits_{k =  - \infty }^{{-1}} {{2^{k(\sum\limits_{i=1}^m\alpha_i(0))p}}\prod\limits_{i = 1}^m {\Big( {\sum\limits_{r = {\Theta _n^*} - 1}^0 {\left\| {{f_i}{\chi _{k + {\ell_i} + r}}} \right\|_{L^{\zeta q_i(\cdot)}_{\omega_i}}} } \Big)^p} } }}} \Big\} ^{\frac{1}{p}}dt.
\nonumber
\end{eqnarray}
We observe that the other estimations can be done by similar arguments as Theorem \ref{TheoremHerzLip}.
Thus, $\mathcal G_0$ and $\mathcal G_1$ are dominated by $\mathcal C.\mathcal C_5.\mathcal B_{CMO,\vec{\omega}}\prod\limits_{i = 1}^m {{{\left\| {{f_i}} \right\|}_{{\mathop K\limits^.}_{\zeta q_i(\cdot),\omega_i}^{{\alpha_i(\cdot)}, p_i}}}}.$
This proves the assertion.
\subsection{Proofs of Theorem \ref{TheoremMorreyLip} and Theorem \ref{TheoremMorreyCMO}}
For $R>0$, we write $B:=B(0,R)$ and $\Delta_R$ as
\[
 \Delta_R={\frac{1}{{{\omega}\big(B\big)^{\frac{1}{q_\infty}+\lambda}}}} \big\|H_{\Phi ,\vec{A}}^{\vec b}(\vec{f})\big\|_{L^{q(\cdot)}_{v}(B)}.
\]
By applying the Minkowski inequality for the variable Lebesgue space, we have
\begin{equation}\label{DeltaR}
\Delta_R \lesssim \int\limits_{\mathbb R^n}{\frac{1}{{\omega(B)}^{\frac{1}{q_\infty} +\lambda}}.\frac{\Phi(t)}{|t|^n}\big\|\prod\limits_{i=1}^{m}{f_i(A_i(t)\cdot)}\big(b_i(\cdot)-b_i(A_i(t)\cdot)\big)\big\|_{L^{q(\cdot)}_{v}(B)}dt}.
\end{equation}
By estimating as (\ref{fLip}) above, we get
\begin{eqnarray}\label{HolderLip}
&&\big\|\prod\limits_{i=1}^{m}{f_i(A_i(t)\cdot)}\big(b_i(\cdot)-b_i(A_i(t)\cdot)\big)\big\|_{L^{q(\cdot)}_{v}(B)}
\nonumber
\\
&&\,\,\,\,\,\,\lesssim R^{\beta}.\mathcal B_{\rm Lip}\prod\limits_{i=1}^m\big\|I_n-A_i(t)\big\|^{\beta_i}.\prod_{i=1}^{m}\big\|f_i(A_i(t)\cdot)\big\|_{L^{q_i(\cdot)}_{v_i}(B)}\big\||\cdot|^{\frac{\alpha_i}{r_i}}\big\|_{L^{r_i}(B)}.
\nonumber
\\
&&\,\,\,\,\,\,\lesssim R^{\beta+\sum\limits_{i=1}^m\frac{\alpha_i+n}{r_i}}.\mathcal B_{\rm Lip}\prod\limits_{i=1}^m\big\|I_n-A_i(t)\big\|^{\beta_i}.\prod_{i=1}^{m}\big\|f_i(A_i(t)\cdot)\big\|_{L^{q_i(\cdot)}_{v_i}(B)}.
\end{eqnarray}
By (\ref{DKnhung1}) and the Theorem \ref{theoembed}, we find
\begin{equation}\label{fiNhung}
\big\|f_i(A_i(t)\cdot)\big\|_{L^{q_i(\cdot)}_{v_i}(B)}\lesssim c_{A_i,q_i,\alpha_i}(t).\big\|1 \big\|_{L^{\vartheta_{1i}(t,\cdot)}}.\big\|f_i\big\|_{L^{q_i(\cdot)}_{v_i}(B(0,R||A_i(t)||))}.
\end{equation}
By the condition (\ref{DKlambda}), we estimate
\begin{equation}\label{trietR}
\frac{R^{\beta+\sum\limits_{i=1}^m\frac{\alpha_i+n}{r_i}}}{{\omega(B)}^{\frac{1}{q_{\infty}} + \lambda }}\lesssim \prod\limits_{i = 1}^m \frac{\big\|A_i(t)\big\|^{(\gamma_i+n)(\frac{1}{q_{i\infty}}+\lambda_i)}}{{\omega_i(B(0,R||A_i(t)||))}^{\frac{1}{{q_{i\infty}}} + {\lambda _i}}}.
\end{equation}
Thus, by having (\ref{DeltaR})-(\ref{trietR}) and defining central Morrey  spaces with variable exponent, it follows that
\[
 \Delta _R \lesssim \mathcal C_6.\mathcal B_{\rm Lip}. \prod\limits_{i = 1}^m {{{\left\| {{f_i}} \right\|}_{{\mathop B\limits^.}_{{\omega_i,v_i}}^{{q_i(\cdot)},{\lambda _i}}}}}
\]
Therefore, we conclude ${\left\| {{H_{\Phi ,\vec{A} }^{\vec b}}(\vec{f} )} \right\|_{{\mathop B\limits^.}_{\omega,v}^{q(\cdot),\lambda }}}\lesssim \mathcal C_{6}.\mathcal B_{\rm Lip}.\prod\limits_{i = 1}^m {{{\left\| {{f_i}} \right\|}_{{\mathop B\limits^.}_{{\omega_i},v_i}^{{q_i(\cdot)},{\lambda _i}}}}}.$
\\

Next, we will prove Theorem \ref{TheoremMorreyCMO}. Indeed, by using the Minkowski inequality and the H\"{o}lder inequality for variable Lebesgue spaces again, it is obvious to show that
\[
\Delta_R \lesssim \int\limits_{\mathbb R^n}{\frac{1}{{\omega(B)}^{\frac{1}{q_\infty} +\lambda}}.\frac{\Phi(t)}{|t|^n}\prod_{i=1}^{m}\big\|b_i(\cdot)-b_i(A_i(t)\cdot)\big\|_{L^{r_i}(v_i,B)}\big\|f_i(A_i(t).)\big\|_{L^{q_i(\cdot)}_{v_i}(B)}dt}.
\]
By (\ref{BDTCMO}) above, we deduce
\begin{eqnarray}
\Delta_R &\lesssim& R^{\sum\limits_{i=1}^m \frac{\alpha_i+n}{r_i}}.\mathcal B_{{CMO,\vec{v}}}\int\limits_{\mathbb R^n}{\frac{1}{{\omega(B)}^{\frac{1}{q_\infty} +\lambda}}.\frac{\Phi(t)}{|t|^n}\prod_{i=1}^{m}\big(1+\psi_{A_i,\alpha_i}^{\frac{1}{r_i}}|s_i(t)|^{\frac{\alpha_i+n}{r_i}}+\varphi_{A_i}(t)\big)}\times
\nonumber
\\
&&\,\,\,\,\,\,\,\,\,\,\,\,\,\,\,\,\,\,\,\,\,\,\,\,\,\,\,\,\,\,\,\,\,\,\,\,\,\,\,\,\,\,\,\,\times\prod\limits_{i=1}^{m}\big\|f_i(A_i(t)\cdot)\big\|_{L^{q_i(\cdot)}_{v_i}(B)}dt.
\nonumber
\end{eqnarray}
For this, by (\ref{fiNhung}), we get
\begin{eqnarray}\label{DeltaCMO}
\Delta_R &\lesssim& R^{\sum\limits_{i=1}^m \frac{\alpha_i+n}{r_i}}\mathcal B_{{CMO,\vec{v}}}
\Big(\int\limits_{\mathbb R^n}{\frac{1}{{\omega(B)}^{\frac{1}{q_\infty} +\lambda}}.\frac{\Phi(t)}{|t|^n}\prod_{i=1}^{m}\big(1+\psi_{A_i,\alpha_i}^{\frac{1}{r_i}}|s_i(t)|^{\frac{\alpha_i+n}{r_i}}+\varphi_{A_i}(t)\big)}\times
\nonumber
\\
&&\times c_{A_i,q_i,\alpha_i}(t)\big\|1 \big\|_{L^{\vartheta_{1i}(t,\cdot)}}.\big\|f_i\big\|_{L^{q_i(\cdot)}_{v_i}(B(0,R||A_i(t)||))}dt\Big).
\end{eqnarray}
On the other hand, by (\ref{DKlambda1}), it follows that
\[
\frac{R^{\sum\limits_{i=1}^m\frac{\alpha_i+n}{r_i}}}{{\omega(B)}^{\frac{1}{q_{\infty}} + \lambda }}\lesssim \prod\limits_{i = 1}^m \frac{\big\|A_i(t)\big\|^{(\gamma_i+n)(\frac{1}{q_{i\infty}}+\lambda_i)}}{{\omega_i(B(0,R||A_i(t)||))}^{\frac{1}{{q_{i\infty}}} + {\lambda _i}}}.
\]
Consequently, by having (\ref{DeltaCMO}), we immediately obtain
$$
{\left\| {{H_{\Phi ,\vec{A} }^{\vec b}}(\vec{f} )} \right\|_{{\mathop B\limits^.}_{\omega, v}^{q(\cdot),\lambda }}}\lesssim \mathcal C_{7}.\mathcal B_{{CMO,\vec{v}}}.\prod\limits_{i = 1}^m {{{\left\| {{f_i}} \right\|}_{{\mathop B\limits^.}_{{\omega_i},v_i}^{{q_i(\cdot)},{\lambda _i}}}}},
$$
which completes the proof.
\\

{\textbf{Acknowledgments}}.  This paper is supported by the Vietnam National Foundation for Science and Technology Development (NAFOSTED) under grant number 101.02-2014.51.

\bibliographystyle{amsplain}

\begin{thebibliography}{79}



\bibitem{Almeida} A. Almeida,  and P. H\"{a}st\"{o}, \textit{Besov spaces with variable smoothness and integrability}, J. Funct. Anal. 258(5) (2010), 1628-1655.

\bibitem{Almeida2012} A. Almeida, D. Drihem, \textit{Maximal, potential and singular type operators on Herz spaces with variable exponents}, J. Math. Anal. Appl. 394 (2012), 781-795.

\bibitem{BM} G. Brown,  F. M\'{o}ricz, \textit{Multivariate Hausdorff operators on the spaces $L^p(\mathbb R^n)$}, J. Math. Anal. Appl. 271 (2002), 443-454.

\bibitem{Bandaliev2010} R.A. Bandaliev, \textit{The boundedness of multidimensional hardy operators in weighted variable Lebesgue spaces}, Lith. Math. J. 50 (2010), 249-259.



\bibitem{Coifman1} R. R. Coifman, Y. Meyer, \textit{On commutators of singular integrals and bilinear singular integrals}, Trans. Amer. Math. Soc. 212 (1975), 315-331.


\bibitem{Chuong} N. M. Chuong, \textit{Degenerate parabolic pseudodifferential operators of variable order}, Dokl. Akad. Nauk SSSR 268 (1983), 1055-1058.

\bibitem{HausdoffCDD} N. M. Chuong, D. V. Duong and K. H. Dung, \textit{Multilinear Hausdorff operators on some function spaces with variable exponent}, arXiv:1709.08185 (2017).


\bibitem{CDH2016} N. M. Chuong, D. V. Duong, H. D. Hung, \textit{Bounds for the weighted Hardy-Ces\`{a}ro operator and its commutator on weighted Morrey-Herz type spaces}, Z. Anal. Anwend. 35 (2016) 489-504.

\bibitem{CHH2016} N. M. Chuong, N. T. Hong, H. D. Hung, \textit{Multilinear Hardy-Cesaro operator and commutator on the product of Morrey-Herz spaces}, Analysis Math. (To appear).

\bibitem{Chuong2016} N. M. Chuong, \textit{Pseudodifferential operators and wavelets over real and p-adic fields}, Springer (Submitted to the Editor of Springer).

\bibitem{CUF2013} D. Cruz-Uribe, A. Fiorenza, \textit{Variable Lebesgue Spaces: Foundations and Harmonic Analysis}, Springer Basel, 2013.

\bibitem{C1982} S. Chanillo, \textit{A note on commutators}, Indiana Univ. Math. J. 31 (1982), 7-16.

\bibitem{CRW1976} R. R. Coifman, R. Rochberg and G. Weiss, \textit{Factorization theorems for Hardy spaces in several
variables}, Ann. of Math. (2), 103 (1976), 611-635.

\bibitem{Diening} L. Diening,  M. R{u}˚\v{z}i\v{c}ka, \textit{Calder\'{o}n-Zygmund operators on generalized Lebesgue spaces $L^{p(x)}$ and problems related to fluid dynamics}, J. Reine Angew. Math. 563 (2003), 197-220. 

\bibitem{Diening1} L. Diening, P. Harjulehto, P. H\"{a}st\"{o}, M. Ruzicka, \textit{Lebesgue and Sobolev spaces with variable exponents}, Springer-Verlag, (2011).

\bibitem{FGLY2015} Z. W. Fu, S. L. Gong, S. Z. Lu and W. Yuan, \textit{Weighted multilinear Hardy operators and commutators}, Forum Math. 27 (2015), 2825-2851.

\bibitem{FLL2009} Z. W. Fu, Z. G. Liu and S. Z. Lu, \textit{Commutators of weighted Hardy operators on $\mathbb R^n$}, Proc Amer. Math. Soc., 137 (2009), pp. 3319-3328.

 \bibitem{Guliyev} V. S. Guliyev, J. Hasanov, and S. Samko, \textit{Boundedness of the maximal, potential and singular operators in the generalized variable exponent Morrey spaces},  Math. Scand., 107 (2010), 285-304.
 
\bibitem{Hausdorff} F. Hausdorff, \textit{Summation methoden und Momentfolgen}, I, Math. Z. 9 (1921), 74-109.

\bibitem{Hurwitz} W. A. Hurwitz, L. L. Silverman, \textit{The consistency and equivalence of certain definitions of summabilities}, Trans. Amer.  Math. Soc. 18 (1917), 1-20.

\bibitem{HK2015} H. D. Hung, L. D. Ky, \textit{New weighted multilinear operators and commutators of Hardy-Ces\`{a}ro type}, Acta Math. Sci. Ser. B Engl. Ed. 35 (2015)(6), 1411-1425.

\bibitem{H2000} W. Hoh, \textit{Pseudodifferential operators with negative definite symbols of varable order}, Revista Mat. Iberoamer. 18, No.2 (2000), 219-241.

\bibitem{I2010} M. Izuki, \textit{Commutators of fractional integrals on Lebesgue and Herz spaces with variable exponent}, Rend. Circ. Mat. Palermo (2) 59 (2010), 461–472.

\bibitem{I2010Hiro} M. Izuki, \textit{Fractional integrals on Herz Morrey spaces with variable exponent}, Hiroshima Math. J. 40 (2010), 343-355.

\bibitem{JN1961} F. John, L. Nirenberg, \textit{On functions of bounded mean oscillation}, Comm Pure Appl. Math., 2 (1961), 415-426.

\bibitem{KL2005} A. Y. Karlovich, A. K. Lerner,  \textit{Commutators of singular integrals on generalized $L ^p$ spaces with variable exponent}, Publ. Mat., 49 (2005), 111-125.

\bibitem{Lu} S. Z. Lu, D. C. Yang, G. E. Hu, \textit{Herz type spaces and their applications}, Beijing Sci. Press (2008).

\bibitem{LY1993} S. Z. Lu, D. C. Yang, \textit{Some new Hardy spaces associated with the Herz spaces and their applications}, J. of Beijing Normal Univ., 29(1993), 10-19.

\bibitem{LZ2014} Y. Lu, Y. P. Zhu, \textit{Boundedness of multilinear Calder\'{o}n-Zygmund singular operators on Morrey-Herz spaces with variable exponents}, Acta Math. Sin.(Engl. Ser.) 30 (2014), 1180-1194.
 

\bibitem{Moricz2005} F. M\'{o}ricz, \textit{Multivariate Hausdorff operators on the spaces $H^1(\mathbb R^n)$ and $BMO(\mathbb R^n)$}, Analysis Math. 31 (2005), 31-41.

\bibitem{Jacob} N. Jacob, H. G. Leopold, \textit{Pseudodifferential operators with variable order of differentiation generating Feller semigroups}, Integr. Equations Oper. Theor. 17 (1993), 544-553.

\bibitem{WZ2016} J. L. Wu, W. J. Zhao, \textit{Boundedness for fractional Hardy-type operator on variable-exponent Herz-Morrey spaces}, Kyoto J. Math. Vol. 56, No. 4 (2016), 831-845.

\bibitem{W2014} J. Wu,  \textit{Boundedness for commutators of fractional integrals on Herz-Morrey spaces with variable exponent}, Kyoto Journal of Mathematics, Vol. 54, No. 3 (2014), 483–495.


\bibitem {Stein} Elias M. Stein, \textit{Harmonic analysis: real-variable methods, orthogonality, and oscillatory integrals,} Princeton University Press, (1993).

\bibitem{TXZ2011} C. Tang, F. Xue and Y. Zhou, \textit{Commutators of weighted Hardy operators on Herz-type spaces}, Ann. Polon. Math. 101 (2011)(3), 267-273.

\bibitem{Xiao2001} J. Xiao, \textit{$L^p$ and $BMO$ bounds of weighted Hardy-Littlewood Averages}, J. Math. Anal. Appl. 262 (2001), 660-666.

\end{thebibliography}

\end{document}